\documentclass[leqno,twoside]{article}
\usepackage{amsmath,amssymb,amsthm}
\usepackage[margin=2.7cm]{geometry} 
\usepackage{titlesec,hyperref}
\usepackage{fancyhdr}

\usepackage{color}

\fancypagestyle{plain}
{
\fancyhf{}

}

\titleformat{\subsection}{\it}{\thesubsection.\enspace}{1pt}{}

\newtheorem{theo}{Theorem}[section]
\newtheorem{lemm}[theo]{Lemma}
\newtheorem{defi}[theo]{Definition}

\newtheorem{prop}[theo]{Proposition}
\newtheorem{rema}[theo]{Remark}
\numberwithin{equation}{section}

\allowdisplaybreaks

\def\th2{\frac{\theta}{2}}

\begin{document}
\title{Well-posedness for the Euler-Nernst-Planck-Possion system\\ in Besov spaces
\hspace{-4mm}
}

\author{Zeng Zhang$^1$
\quad Zhaoyang Yin$^2$ \\[10pt]
Department of Mathematics, Sun Yat-sen University,\\
510275, Guangzhou, P. R. China.\\[5pt]
}
\footnotetext[1]{Email: \it zhangzeng534534@163.com}
\footnotetext[2]{Email: \it mcsyzy@mail.sysu.com.cn}
\date{}
\maketitle

\begin{abstract}
In this paper, we mainly study the Cauchy problem of the Euler-Nernst-Planck-Possion $(ENPP)$ system. We first establish local well-posedness for the Cauchy problem of the $ENPP$ system in Besov spaces. Then we present a blow-up criterion of solutions to the $ENPP$ system. Moreover, we prove that the solutions of the Navier-Stokes-Nernst-Planck-Possion system converge to the solutions of the  $ENPP$ system as the viscosity $\nu$ goes to zero, and that the convergence rate is at least of order ${\nu}^\frac{1}{2}$.

\vspace*{5pt}

\vspace*{5pt}
\end{abstract}

\vspace*{10pt}
\tableofcontents

\section{Introduction}
Electro-kinetics describes the dynamic coupling between incompressible flows and diffuse charge systems. It finds application in biology, chemistry and pharmacology \cite{Joseph}. We shall study such a model in this paper
 as follows:
\begin{align}
u_t+u\cdot \nabla u-\nu \triangle u+\nabla P&=\triangle\phi\nabla\phi, \quad t>0,\,x \in \mathbb{R}^d, \label{1.1}\\[1ex]
\nabla\cdot u&=0,\quad t>0,\,x \in \mathbb{R}^d,\\[1ex]
n_t+u\cdot \nabla n&=\nabla\cdot(\nabla n-n\nabla\phi), \quad t>0,\,x \in \mathbb{R}^d,\label{1.3}\\[1ex]
p_t+u\cdot \nabla p&=\nabla\cdot(\nabla p+p\nabla\phi), \quad t>0,\,x \in \mathbb{R}^d,\label{1.4}\\[1ex]
\triangle\phi&=n-p,\quad t>0,\,x \in \mathbb{R}^d,\label{1.5}\\[1ex]
(u,n,p)|_{t=0}&=(u_0,n_0,p_0),\quad x \in \mathbb{R}^d.
\end{align}
Here $u(t,x)$ is the velocity of the fluid, $P$ is the pressure, $\nu\geq 0$ is the fluid viscosity, $n$ and $p$ are the densities of the negative and positive charged particles, and $\phi$ is the electronic potential. The first two equations represent the momentum equation and the incompressibility
of the solution, and the right hand side of (\ref{1.1}) is the Lorentz (or Coulomb) force caused by the charges. (\ref{1.3}) and (\ref{1.4}) model the balance between diffusion and convective transport of charge densities by flow and electric fields. (\ref{1.5}) is the Possion equation for the electrostatic potential $\phi$ where the right hand side is the net charge density. The above system (1.1)-(1.6) is the so called Navier-Stokes-Nernst-Planck-Possion ($NSNPP$) system if $\nu>0,$ and the Euler-Nernst-Planck-Possion ($ENPP$) system if $\nu=0$ \cite{Joseph}.

It is obvious that if the flow is charge free, i.e. $n=p=0$, the system (1.1)-(1.6) reduces to the Navier-Stokes/Euler system. If, on the other hand, in the absence of a fluid, a.e. $u=0$,  the system (1.1)-(1.6) reduces to the Nernst-Planck-Possion system, which was first studied by Nernst and Plank at the end of the nineteenth century as a basic model for the diffusion of ions in an electrolyte filling all of $\mathbb{R}^3,$ see \cite{Biler,H, Selberherr} for more details.

There are several mathematical results on the $NSNPP$ system. By using Kato's semigroup framework, Joseph \cite{Joseph} obtained the existence of a unique smooth local solution for smooth initial dada, and established the stability under the inviscid limit to the $ENPP$ system. Applying Schauder's fixed point theorem,
Schmuck \cite{Schmuck} established the global existence of weak solutions and local existence of strong solutions in a bounded domain under some boundary and initial conditions. Ryham \cite{Ryham} proved the existence, uniqueness, and regularity of weak solutions in a bounded domain. Zhao et al. \cite{zhaoW,zhaoW3,zhaoG,zhaoW2} studied well-posedness in Lebesgue spaces, modulation spaces, Triebel-Lizorkin spaces and Besov spaces by using the Banach fixed point theorem. Li \cite{Li} studied the quasineutral limit by establishing elaborate energy analysis.

However, to the best of the authors' knowledge, there are no mathematical discussions on the $ENPP$ system.
In this paper, we study well-posedness for the Cauchy problem of the $ENPP$ system in Besov spaces in dimension $d\geq2.$ The main difficulty here comes from the term $\nabla\phi$ determined by the Possion equation (\ref{1.5}). In fact, if $n-p\in L^{p}$ with $1<p<d,$
 due to the Hardy-Littlewood-Sobolev inequality, $\nabla\phi=\nabla(-\Delta)^{-1}(p-n)\in L^{\frac{pd}{p+d}}.$ We will establish more rigorous product laws in Besov spaces, to make the product terms $\triangle\phi\nabla\phi,$ $\nabla\cdot(n\nabla\phi)$ and $\nabla\cdot(p\nabla\phi)$ have the same regularities as the other terms in the associated equations.

Before stating our main results, we define spaces $L^\infty_L(\mathbb{R}^d),~L^\infty_{P^\alpha}(\mathbb{R}^d)(0<\alpha<1)$ and their norms by
\begin{align*}
  L^\infty_L(\mathbb{R}^d)=\{u\big|\|u\|_{L^\infty_L}=\sup_{x\in \mathbb{R}^d}\frac{|u(x)|}{1+log\langle x\rangle}<\infty\},\\
  L^\infty_{P^\alpha}(\mathbb{R}^d)=\{u\big|\|u\|_{L^\infty_{P^\alpha}}=\sup_{x\in \mathbb{R}^d}\frac{|u(x)|}{\langle x\rangle^\alpha}<\infty\},
\end{align*}
where $\langle x\rangle=(1+|x|^2)^{\frac{1}{2}}.$
It is obvious that $L^\infty_L(\mathbb{R}^d)~and~L^\infty_{P^\alpha}(\mathbb{R}^d)(0<\alpha<1)$ are Banach spaces.\\
We also denote $X(T)$ the set of functions $(a,b,c)$ in
\begin{align*}
\widetilde{L}^\infty_{T}(B^{s_1}_{p_1,r_1}(\mathbb{R}^d))\times \Big(\widetilde{L}^\infty_{T}(B^{s_2}_{p_2,r_2}(\mathbb{R}^d))
\cap\widetilde{L}^1_{T}(B^{s_2+2}_{p_2,r_2}(\mathbb{R}^d))\Big)^2
\end{align*}
endowed with the norm
\begin{align*}
\|(a,b,c)\|_{X(T)}=\|a\|_{\widetilde{L}^\infty_{T}(B^{s_1}_{p_1,r_1}(\mathbb{R}^d))}
+\|b\|_{\widetilde{L}^\infty_{T}(B^{s_2}_{p_2,r_2}(\mathbb{R}^d))
\cap\widetilde{L}^1_{T}(B^{s_2+2}_{p_2,r_2}(\mathbb{R}^d))}+\|c\|_{\widetilde{L}^\infty_{T}(B^{s_2}_{p_2,r_2}(\mathbb{R}^d))
\cap\widetilde{L}^1_{T}(B^{s_2+2}_{p_2,r_2}(\mathbb{R}^d))},
\end{align*}
and $Y_{\alpha}(T)$ the set of functions $(e,f)$ in
\begin{align*}
L^1_{T}(L^\infty_L(\mathbb{R}^d))\times L^\infty_T(L^\infty_{P^\alpha}(\mathbb{R}^d))
\end{align*}
endowed with the norm
\begin{align*}
\|(e,f)\|_{Y_{\alpha}(T)}=\|e\|_{L^1_{T}(L^\infty_L(\mathbb{R}^d))}
+\|f\|_{L^\infty_T(L^\infty_{P^\alpha}(\mathbb{R}^d))}.
\end{align*}

We can now state our main results:

\begin{theo}\label{a1}
Let $d\geq 2,~(s_1,s_2)\in\mathbb{R}^2,~~1\leq r_1,r_2\leq\infty,~1<p_1\leq\infty,~and~1<p_2<d,$ satisfying
\begin{align}
  s_1&>1+\frac{d}{p_1},&\textit{or}~s_1=&1+\frac{d}{p_1}~\textit{and}~r=1,&\textit{and}~&(s_1,p_1,r_1)\neq(1,\infty,1),\label{1.7}\\
  s_2&<s_1<s_2+\frac{3}{2},&s_2-\frac{d}{p_2}<&s_1-\frac{d}{p_1}<s_2-\frac{d}{p_2}+\frac{3}{2},\label{1.8}\\
  \frac{1}{p_1}&\leq\frac{2}{p_2}-
   \frac{1}{d},&\frac{1}{r_1}\leq&\frac{2}{r_2}.\label{1.9}
\end{align}
There exists constants $c$ and $r\geq 4$, depending only on $s_1,p_1,r_1,s_2,p_2,r_2$ and $d,$ such that for $(u_0,n_0,p_0)\in B^{s_1}_{p_1,r_1}(\mathbb{R}^d)\times (B^{s_2}_{p_2,r_2}(\mathbb{R}^d))^2,$ and ${\rm div}u_0=0,~n_0\geq 0, p_0\geq 0,$ there exists a time $$T\geq \frac{c}{1+(\|u_0\|_{ B^{s_1}_{p_1,r_1}(\mathbb{R}^d)}+\|n_0\|_{B^{s_2}_{p_2,r_2}(\mathbb{R}^d)}+\|p_0\|_{B^{s_2}_{p_2,r_2}(\mathbb{R}^d)})^r},$$ such that the $ENPP$ system has a solution $(u,n,p,P,\phi)$ on $[0,T]\times\mathbb{R}^d$ satisfying
$$(u,n,p)\in X(T),~(P,\phi)\in Y_\alpha(T)~ \textit{for some}~\alpha\in(0,1).$$
Moreover, $n,p\geq 0,~a.e.~on ~[0,T]\times \mathbb{R}^d,$ and \\ $u\in C([0,T];B^{s_1}_{p_1,r_1}(\mathbb{R}^d)),$ if $r_1<\infty,$ or $u\in C([0,T];B^{\tilde{s}_1}_{p_1,r_1}(\mathbb{R}^d)),$ if $r_1=\infty,~\tilde{s}_1<s_1,$\\ $(n,p)\in \big(C([0,T];B^{s_2}_{p_2,r_2}(\mathbb{R}^d))\big)^2,$ if $r_2<\infty,$ or $(n,p)\in  \big(C([0,T];B^{\tilde{s}_2}_{p_2,r_2}(\mathbb{R}^d))\big)^2,$ if $r_2=\infty,~\tilde{s}_2<s_2.$\\
Finally, if $(\tilde{u},\tilde{n},\tilde{p},\tilde{P},\tilde{\phi})$ also solves the $ENPP$ system with the same initial data and belongs to $X(T)\times Y_{\alpha_1}(T)$ for some $\alpha_1\in(0,1)$, then
$$\tilde{u}=u,~\tilde{n}=n,~\tilde{p}=p,\nabla\tilde{P}=\nabla P,~\nabla\tilde{\phi}=\nabla \phi.$$
\end{theo}
\begin{rema}
We mention that the restriction $(\ref{1.7})$ is due to some reasons as illustrated for the Euler equation in \cite{keben}, and the conditions $(\ref{1.8})$ and $(\ref{1.9})$ are caused by the coupling between $u$ and $(n,p)$ and the product laws in Besov spaces.
We also point out that from the condition $(\ref{1.8})$, there exists some $\varepsilon_1\in(0,\frac{3}{8})$ such that
\begin{align}\label{ee}
s_2+4\varepsilon_1<s_1+2\varepsilon_1<s_2+\frac{3}{2},
~~~s_2+4\varepsilon_1-\frac{d}{p_2}<s_1+2\varepsilon_1-\frac{d}{p_1}<s_2-\frac{d}{p_2}+\frac{3}{2}.
\end{align}
\end{rema}
\begin{rema}
Note that for every $s\in \mathbb{R},$ $H^s(\mathbb{R}^d)=B^s_{2,2}(\mathbb{R}^d).$ For $d\geq 3,$ Theorem \ref{a1} holds true in Sobolev spaces $H^{s_1}(\mathbb{R}^d)\times \big(H^{s_2}(\mathbb{R}^d)\big)^2$ with $s_1>1+\frac{d}{2},~s_2<s_1<s_2+\frac{3}{2}.$
\end{rema}
\begin{theo}\label{a2}
Under the assumptions of Theorem \ref{a1}, assume further
 $s_2>\frac{d}{p_2},$ or $s_2=\frac{d}{p_2},~r_2=1,$ and that the $ENPP$ system has a solution $(u,n,p,P,\phi)\in X(T)\times Y_\alpha(T)$ for some $\alpha\in (0,1).$  If
\begin{align}\label{b1}
\left\{
\begin{array}{ll}
\|u\|_{L^\infty_T(L^\infty(\mathbb{R}^d))}+\int_0^T \|\nabla u(t)\|_{L^\infty(\mathbb{R}^d)}dt
<\infty,~~~p_1=\infty,\\ [1ex]
\int_0^T\|\nabla u(t)\|_{L^\infty(\mathbb{R}^d)}dt
<\infty,~~~1<p_1<\infty,\end{array}\right.
\end{align}
then there exists some $T^*>T,$ such that $(u,n,p,P,\phi)$ can be continued on $[0,T^*]\times \mathbb{R}^d$ to a solution of the $ENPP$ system which belongs to $X(T^*)\times Y_\alpha(T^*).$
\end{theo}
\begin{rema}
 For $d\geq 3,$ Theorem \ref{a2} holds true in Sobolev spaces $H^{s_1}(\mathbb{R}^d)\times \big(H^{s_2}(\mathbb{R}^d)\big)^2$ with $s_1>1+\frac{d}{2},~s_2>\frac{d}{2},~s_2<s_1<s_2+\frac{3}{2}.$
\end{rema}
\begin{theo}\label{a3}
 Under the assumptions of Theorem \ref{a1}, there exist positive constants $T$ and $M$ independent of $\nu,$ such that the
  $NSNPP$ system
  has a unique solution  $(u_\nu,n_\nu,p_\mu,P_\nu,\Psi_\nu)$ in $X(T)\times Y_\alpha(T)$ for some $\alpha\in(0,1).$ Moreover, the solution $(u_\nu,n_\nu,p_\mu,P_\nu,\Psi_\nu)$ converges to a solution $(u,n,p,P,\Psi)$ to the $ENPP$ system as the viscosity $\nu$ goes to zero, and the convergence rate of $$\|(u_{\nu}-u,n_{\nu}-n,p_{\nu}-p)\|_{\widetilde{L}^\infty(B^{s'_1}_{p_1,r_1})\times \big(\widetilde{L}^\infty(B^{s_2-1}_{p_2,r_2})\big)^2}$$ is at least of order ${\nu}^\frac{1}{2},$ where $s'_1$ satisfy
  \begin{align}\label{s'}
\left\{\begin{array}{l}
s_1'=s_1-1-\varepsilon_1,~~if~s_1=2+\frac{d}{p}, \\
s_1'=s_1-1,~~~~~~~~~otherwise,
\end{array}\right.
\end{align}
with $\varepsilon_1$ defined as in (\ref{ee}).
\end{theo}
Throughout the paper, $C>0$ denotes various ``harmless" finite constant, $c>0$ denotes a small constant. We shall sometimes use $X\lesssim Y$ to denote $X\leq CY.$ For simplicity, we write $L^P,~L^\infty_L,~L^\infty_{P^\alpha}$ and $B^s_{p,r}$ for the spaces $L^P(\mathbb{R}^d),~L^\infty_L(\mathbb{R}^d),~L^\infty_{P^\alpha}(\mathbb{R}^d)$ and $B^s_{p,r}(\mathbb{R}^d)$, respectively. We mention that according to the context, $p$ denotes as the index of the Besov space or the density of the positive charged particles is not confused.

  The remain part of this paper is organized as follows. We introduce Besov spaces and the modified $ENPP$ system in Section 2. In Section 3, we prove Theorem \ref{a1} by using a more accurate product estimate. Section 4 is devoted to the proof of Theorem \ref{a2}. Finally, in Section 5, we prove Theorem \ref{a3}.
\section{Preliminaries}
  \subsection{The nonhomogeneous Besov spaces }
We first define the Littlewood-Paley decomposition.
 \begin{lemm}\cite{keben}
Let $\mathcal{C}=\{\xi\in{\mathbb{R}^2}, ~\frac{3}{4}\leq|\xi|\leq\frac{8}{3} \}$ be an annulus. There exist radial functions $\chi$ and $\varphi$ valued in the interval $[0,1]$, belonging respectively to $\mathcal{D}(B(0,\frac{4}{3}))$ and $\mathcal{D}(\mathcal{C})$,  such that
\begin{align*}
\forall \xi \in{\mathbb{R}}^d,~~\chi(\xi)+\sum_{j\geq 0}\varphi(2^{-j}\xi)=1.\end{align*}\end{lemm}
The nonhomogeneous dyadic blocks $\triangle_j$ and the nonhomogeneous low-frequency cut-off operator $S_j$ are then defined as follows:
\begin{align*}&\triangle_ju=0~~if~j\leq-2,~~~~~~~~~~~~~\triangle_{-1}u=\chi(D)u,\\
&\triangle_ju=\varphi(2^{-j}D)u~~if~j\geq0,~~~~S_ju=\sum_{j'\leq j-1}\triangle_{j'}u,~~for~j\in \mathbb{Z}.\end{align*}

We may now introduce the nonhomogeneous Besov spaces.
\begin{defi}\label{dingyi}
Let $s\in \mathbb{R}$ and $(p,r)\in[1,\infty]^2$. The nonhomogeneous Besov space $B^s_{p,r}$
consists of all tempered distributions $u$ such that
$$\|u\|_{B^s_{p,r}}\overset{def}{=}\Big\|(2^{js}\|\triangle_ju\|_{L^p})_{j\in \mathbb{Z}}\Big\|_{l^r(\mathbb{Z})}
<\infty.$$
\end{defi}
\begin{lemm}\label{Fadou}
The set $B^s_{p,r}$ is a Banach space, and satisfies the Fatou property, namely, if $(u_n)_{n\in N}$ is a
bounded sequence of $B^s_{p,r}$, then an element $u$ of $B^s_{p,r}$ and a subsequence $u_{\psi(n)}$ exist such that
$$\underset{n\rightarrow\infty}{\lim}~u_{\psi(n)}=u~~in~~\mathcal{S}'~~~and ~~~\|u\|_{B^s_{p,r}}\leq C \underset{n\rightarrow\infty}{\liminf} \|u_{\psi(n)}\|_{B^s_{p,r}}.$$\end{lemm}
In addition to the general time-space $L^{\rho}_T(B^s_{p,r})$, we introduce the following mixed time-space $\widetilde{L}^{\rho}_T(B^s_{p,r}).$
\begin{defi}
For all $T>0,~s\in\mathbb{R},$ and $1\leq r,\rho\leq\infty$, we define the space $\widetilde{L}^{\rho}_T(B^s_{p,r})$ the set of tempered distributions $u$ over $(0,T)\times \mathbb{R}^d,$ such that
$$\|u\|_{\widetilde{L}^{\rho}_T(B^s_{p,r})}\overset{def}{=}\|2^{js}\|\triangle_ju\|_{L^{\rho}_T(L^p)}\|
_{l^r(\mathbb{Z})}<\infty.$$
\end{defi}
\noindent It follows from the Minkowski inequality that
\begin{align*}
  \|u\|_{L^\rho_T({B}_{p,r}^{s})}\leq \|u\|_{\widetilde{L}^\rho_T({B}_{p,r}^{s})}~if~r\leq \rho,~~\|u\|_{\widetilde{L}_T^\rho({B}_{p,r}^{s})}\leq \|u\|_{L^\rho_T({B}_{p,r}^{s})}~if~r\geq \rho.
\end{align*}

  Let's then recall Bernstein-Type lemmas.
\begin{lemm}\label{Bi}\cite{keben} (Bernstein inequalities)
Let $\mathcal{C}$ be an annulus and $\mathcal{B}$ a ball. A constant $C$ exists such that for any nonnegative integer $k$, any couple $(p,q)$ in $[1,\infty]^2$ with $q\geq p\geq1$, and any function u of $L^p$, we have
\begin{align*}
&Supp\, \widehat{u}\subset \lambda\mathcal{B}\Rightarrow\,\underset{|\alpha|=k}{\sup}\,\|\partial^{\alpha}u\|_{L^q}
\leq C^{k+1}\lambda^{k+d(\frac{1}{p}-\frac{1}{q})}\|u\|_{L^p},\\
&Supp\,\widehat{u}\subset \lambda\mathcal{C}\Rightarrow\,C^{-k-1}\lambda^k\|u\|_{L^p}\leq\underset{|\alpha|=k}{\sup}\,\|\partial^{\alpha}u\|_{L^q}
\leq C^{k+1}\lambda^k\|u\|_{L^p}.\end{align*}
\end{lemm}
\begin{lemm}\label{DL}\cite{DanchinL} (A Bernstein-like inequality)
Let $1< p <\infty$ and $u\in L^p(\mathbb{R}^d),$ such that $Supp
\,
\widehat{u}\in C(0,R_1,R_2)~(with~ 0 < R_1 < R_2).$ There exists a constant c depending
only on $d$ and $R_2/R_1$, such that
\begin{align*}
c\frac{R_1^2}{p^2}\int_{\mathbb{R}^d}|u|^pdx\leq\int_{\mathbb{R}^d} |\nabla u|^2|u|^{p-2}dx=-\frac{1}{p-1}\int_{\mathbb{R}^d}
\triangle u|u|^{p-2}udx.
\end{align*}
\end{lemm}

 We state the following embedding and interpolation inequalities.
\begin{lemm}\cite{keben}
  Let $1\leq p_1\leq p_2\leq\infty$ and $\leq r_1\leq r_2\leq\infty.$ Then for any real number $s,$ we have
  $${B}^s_{p_1,r_1}\hookrightarrow{B}^{s-d(\frac{1}{p_1}-\frac{1}{p_2})}_{p_2,r_2}.$$
  \end{lemm}
  \begin{lemm}\cite{keben}
  If $s_1$ and $s_2$ are real numbers such that $s_1<s_2,$  $\theta\in(0,1)$ and $~1\leq p,r\leq\infty,$ then we have
 \begin{align*}
 B^{s_2}_{p,\infty}\hookrightarrow B^{s_1}_{p,1},~~~\textit{and}~~~
 \|u\|_{{B}^{\theta s_1+(1-\theta)s_2}_{p,r}}\leq \|u\|_{B^{s_1}_{p,r}}^\theta\|u\|_{B^{s_2}_{p,r}}^{1-\theta}.
  \end{align*}
  \end{lemm}

In the sequel, we will frequently use the Bony decomposition:
$$uv=T_vu+T_uv+R(u,v),$$
with
\begin{align*}&R(u,v)=\underset{|k-j|\leq1}{\sum}\triangle_ku\triangle_jv,\\
&T_uv=\underset{j\in \mathbb{Z}}{\sum}S_{j-1}u\triangle_jv
=\underset{j\geq1}{\sum}S_{j-1}u\triangle_j\big((Id-\triangle_{-1})v\big)~~
,\end{align*}
where operator $T$ is called ``paraproduct", whereas $R$ is called ``remainder".
\begin{lemm}\label{T}
A constant $C$ exists which satisfies the following inequalities for any couple of real numbers $(s,t)$ with t negative and any $(p,p_1,p_2,r,r_1,r_2)$ in $[1,\infty]^6$:
\begin{align*}&\|T\|_{\mathcal{L}(L^{p_1}\times {B}^s_{p_2,r};{B}^s_{p,r})}\leq C^{|s|+1},\\
&\|T\|_{\mathcal{L}({B}^t_{p_1,r_1}\times {B}^s_{p_2,r_2};{B}^{s+t}_{p,r})}\leq \frac{C^{|s+t|+1}}{-t},
\end{align*}
with $\frac{1}{p}\overset{def}{=}\frac{1}{p_1}+\frac{1}{p_2}\leq1,~ \frac{1}{r}\overset{def}{=}min\{1,\frac{1}{r_1}+\frac{1}{r_2}\}.$
\end{lemm}
\noindent{Proof.} The proof of this lemma can be easily deduced from substituting the estimate
$$\|S_{j-1}u\triangle_j v\|_{L^p}\leq \|S_{j-1}u\|_{L^{p_1}}\|\triangle_j v\|_{L^{p_2}},$$
for the estimate
$$\|S_{j-1}u\triangle_j v\|_{L^p}\leq \|S_{j-1}u\|_{L^{\infty}}\|\triangle_j v\|_{L^{p}}$$
in the proof of Theorem 2.82 in \cite{keben}. It is thus omitted.\qed
\begin{lemm}\label{R}\cite{keben}
A constant $C$ exists which satisfies the following inequalities. Let $(s_1,s_2)$ be in $\mathbb{R}^2$ and
$(p_1,p_2,r_1,r_2)$ be in $[1,\infty]^4$. Assume that
$$\frac{1}{p}\overset{def}{=}\frac{1}{p_1}+\frac{1}{p_2}\leq1~~and~~
\frac{1}{r}\overset{def}{=}\frac{1}{r_1}+\frac{1}{r_2}\leq1.$$
If $s_1+s_2>0$, then we have, for any $(u,v)$ in ${B}^{s_1}_{p_1,r_1}\times {B}^{s_2}_{p_2,r_2}$,
$$\|R(u,v)\|_{{B}^{s_1+s_2}_{p,r}}\leq\frac{C^{|s_1+s_2|+1}}{s_1+s_2}\|u\|_{{B}^{s_1}_{p_1,r_1}}\|v\|_{{B}^{s_2}_{p_2,r_2}}.$$
If $r=1$ and $s_1+s_2=0$, then we have, for any $(u,v)$ in ${B}^{s_1}_{p_1,r_1}\times {B}^{s_2}_{p_2,r_2}$,
$$\|R(u,v)\|_{{B}^0_{p,\infty}}\leq C\|u\|_{{B}^{s_1}_{p_1,r_1}}\|v\|_{{B}^{s_2}_{p_2,r_2}}.$$
\end{lemm}
We mention that all the properties of continuity for the paraproduct and remainder remain true in the mixed time-space $\widetilde{L}^{\rho}_T(B^s_{p,r}).$

Finally, we state the following commutator estimates.
\begin{lemm}\label{jiaohuan}\cite{keben}
 Let $v$ be a vector filed over $\mathbb{R}^d,$ define $R_j=[v\cdot\nabla,\triangle_j]f.$ Let $\sigma>0~(\textit{or}~\sigma>-1,~if~{\rm div}~v=0),$ $1\leq r\leq\infty,$ $1\leq p\leq p_1\leq\infty,$ and $\frac{1}{p_2}=\frac{1}{p}-\frac{1}{p_1}.$ Then
 \begin{align*}
   \Big\|2^{j\sigma}\|R_j\|_{L^P}\Big\|_{l^r}\leq C\Big(\|\nabla v\|_{L^\infty}\|f\|_{B^{\sigma}_{p,r}}+\|\nabla f\|_{L^{p_2}}\|\nabla v\|_{B^{\sigma-1}_{p_1,r}}\Big).
 \end{align*}
\end{lemm}
\subsection{A priori estimates for transport and transport-diffusion equations}
Let us state some classical a priori estimates for transport equations and transport-diffusion equations.
\begin{lemm}\label{ts}\cite{keben}
Let $1\leq p\leq p_1\leq\infty,~1\leq r\leq\infty$. Assume that
\begin{align}\label{tiaojian}
s\geq-d\,min\left(\frac{1}{p_1},\frac{1}{p'}\right) \quad \textit{or} \quad s\geq-1-d\,min\left(\frac{1}{p_1},\frac{1}{p'}\right)~\textit{if}~div~ v=0
\end{align}
with strict inequality if $r<\infty$.

There exists a constant $C$, depending only on $d, p, p_1, r$ and $s$, such that for all solutions
$f\in L^{\infty}([0,T];B^s_{p,r})$ of the transport equation
\begin{align}
\left\{
\begin{array}{l}
\partial_tf+v\cdot\nabla f=g\\
f_{|t=0}=f_0,
\end{array}
\right.
\end{align}
with initial data $f_0$ in $B^s_{p,r}$, and $g$ in $L^1([0,T];B^s_{p,r})$, we have, for $a.e.\,t\in[0,T]$,
\begin{align}\label{,}\|f\|_{\widetilde{L}_t^{\infty}(B^s_{p,r})}\leq\left(\|f_0\|_{B^s_{p,r}}+
\int_0^t exp(-CV_{p_1}(t))\|g(t')\|_{B^s_{p,r}}dt'\right)exp(CV_{p_1}(t)),\end{align}
or more accurately,
\begin{align}\label{,,}\|f\|_{\widetilde{L}_t^{\infty}(B^s_{p,r})}\leq\|f_0\|_{B^s_{p,r}}+
C\int_0^t V_{p_1}(t')\|f(t')\|_{\widetilde{L}_{t'}^{\infty}(B^s_{p,r})}dt'+
\|g\|_{\widetilde{L}_t^{1}(B^s_{p,r})},\end{align}
with, if the inequality is strict in (\ref{tiaojian}),
\begin{align}
V'_{p_1}(t)=\left\{\begin{array}{l}\|\nabla v(t)\|_{B^{s-1}_{p_1,r}},~if~s>1+\frac{d}{p_1}~or~s=1+\frac{d}{p_1},~r=1,\\
\|\nabla v(t)\|_{B^{\frac{d}{p_1}}_{p_1,\infty}\cap L^{\infty}},~if~s<1+\frac{d}{p_1}
\end{array}\right.
\end{align}
and, if equality holds in (\ref{tiaojian}) and $r=\infty$,
$$V'_{p_1}=\|\nabla v(t)\|_{B^{\frac{d}{p_1}}_{p_1,1}}.$$
If $f=v$, then for all $s>0$ $(s>-1,$ if $div\,u=0)$, the estimates (\ref{,}) and (\ref{,,}) hold with
$$V'_{p_1}(t)=\|\nabla u\|_{L^{\infty}}.$$
\end{lemm}
\begin{lemm}\label{dts}\cite{keben}
Let $1\leq p_1\leq p\leq\infty,~1\leq r\leq\infty,~s\in\mathbb{R}$ satisfy (2.10), and let $V_{p_1}$ be defined as in Lemma \ref{ts}.

There exists a constant $C$ which depends only on $d, r, s$ and $s-1-\frac{d}{p_1}$ and is such that for any smooth solution $f$ of the transport diffusion equation
\begin{align}
\left\{
\begin{array}{l}
\partial_tf+v\cdot\nabla f-\nu\triangle f=g\\
f_{|t=0}=f_0,
\end{array}
\right.
\end{align} we have
\begin{align*}\nu^{\frac{1}{\rho}}\|f\|_{\widetilde{L}^{\rho}_T(B^{s+\frac{2}{\rho}}_{p,r})}\leq Ce^{C(1+\nu T)^{\frac{1}{\rho}}V_{p_1}(T)}\Big(&(1+\nu T)^{\frac{1}{\rho}}\|f_0\|_{B^s_{p,r}}\\
+&(1+\nu T)^{1+\frac{1}{\rho}-\frac{1}{\rho_1}}\nu^{\frac{1}{\rho_1}-1}\|g\|_
{\widetilde{L}^{\rho_1}_T(B^{s-2+\frac{2}{\rho_1}}_{p,r})}\Big),
\end{align*} where $1\leq\rho_1\leq\rho\leq\infty.$
\end{lemm}
\subsection{The modified $ENPP$ system}
Motivated by \cite{keben} for the study of the Euler system, we also introduce the following modified ENPP system
\begin{align}\tag{$\widetilde{ENPP}$}
\left\{
\begin{array}{c}
u_t+u\cdot \nabla u+\Pi(u,u)=\mathcal{P}\big((n-p)\nabla(-\triangle)^{-1}(p-n)\big),  \\[1ex]
n_t-\triangle n=-\big(T_u\nabla n+T_{\nabla n}u+\nabla\cdot R(un))-\nabla\cdot(n\nabla(-\triangle)^{-1}(p-n)\big), \\[1ex]
p_t-\triangle p=-\big(T_u\nabla p+T_{\nabla p}u+\nabla\cdot R(up))+\nabla\cdot(p\nabla(-\triangle)^{-1}(p-n)\big),\\[1ex]
(u,n,p)|_{t=0}=(u_0,n_0,p_0),
\end{array}
\right.
\end{align}
where $\mathcal{P}$ is the Leray projector defined as $\mathcal{P}=Id+\nabla (-\triangle)^{-1} div,$ and $\Pi(\cdot,\cdot)$ is a bilinear operator defined by
\begin{align*}
  \Pi(u,v)=\sum_{j=1}^{5}\Pi_j(u,v),
\end{align*}
where $$\Pi_1(u,v)=\nabla|D|^{-2}T_{\partial_iu^j} \partial_jv^i,$$
$$\Pi_2(u,v)=\nabla|D|^{-2}T_{\partial_jv^i} \partial_iu^j,$$
$$\Pi_3(u,v)=\nabla|D|^{-2}\partial_i\partial_j(I-\triangle_{-1})R(u^i,v^j),$$
$$\Pi_4(u,v)=\theta E_d\ast\nabla \partial_i\partial_j\triangle_{-1}R(u^i,v^j),$$
$$\Pi_5(u,v)=\nabla \partial_i\partial_j\big((1-\theta ) E_d\big)\ast\triangle_{-1}R(u^i,v^j),$$
here $\theta$ is a function of $\mathcal{D}(B(0,2))$ with value $1$ on $B(0,1),$ $E_d$ stands for the fundamental solution of $-\triangle,$ and $|D|^{-2}$ denotes the Fourier multiplier with symbol $|\xi|^2.$

We recall some basic results for $\Pi(\cdot,\cdot).$ See \cite{keben} (Pages 296-300) for further details.
\begin{lemm}\label{pi}\cite{keben}
  For all $s>-1,$ and $1\leq p,r\leq\infty,$ there exists a constant C such that
  \begin{align*}
    \|\Pi(v,w)\|_{B^s_{p,r}}\leq C(\|v\|_{C^{0,1}}\|w\|_{B^s_{p,r}}+\|w\|_{C^{0,1}}\|v\|_{B^s_{p,r}}).
  \end{align*}
  Moveover, there exists a bilinear operator $P_\Pi$ such that
   $\Pi(v,w)=P_\Pi(v,w),$ and
   \begin{align*}
    \|P_\Pi(v,w)\|_{B^{s+1}_{p,r}}\leq C\Big(\|v\|_{C^{0,1}}\|w\|_{B^s_{p,r}}+\|v\|_{C^{0,1}}\|v\|_{B^s_{p,r}}\Big),~if ~ 1<p<\infty,\\
    \|P_\Pi(v,w)\|_{B^{s+1}_{p,r}+L^\infty_L}\leq C\Big(\|v\|_{C^{0,1}}\|w\|_{B^s_{p,r}}+\|v\|_{C^{0,1}}\|v\|_{B^s_{p,r}}\Big),~if ~ p=\infty.
  \end{align*}
\end{lemm}
\begin{lemm}\label{deng}\cite{keben}
  For all $-1<s<\frac{d}{p}+1,$ and $1\leq p,r\leq\infty,$ we have
  \begin{align*}
    \|\Pi(v,w)\|_{B^s_{p,r}}\leq C\Big(\|v\|_{C^{0,1}}\|w\|_{B^s_{p,r}}+\|w\|_{B^{s-\frac{d}{p}}_{\infty,\infty}}\|\nabla v\|_
    {B^{\frac{d}{p}}_{p,r}}\Big).
  \end{align*}
  \end{lemm}
  \begin{lemm}\label{yyy}\cite{keben}
  For all $s>1,$ and $1\leq p,r\leq\infty,$ there exists a constant C such that
  \begin{align*}
    \|div~\Pi(v,w)+tr(Dv,Dw)\|_{B^{s-1}_{p,r}}\leq C\Big(\|div~v\|_{B^{0}_{\infty,\infty}}\|w\|_{B^s_{p,r}}+
    \|div~w\|_{B^{0}_{\infty,\infty}}\|v\|_{B^s_{p,r}}\Big).
  \end{align*}
  \end{lemm}
 \begin{lemm} \label{bao}\cite{keben}
  Let $1 < p < \infty$. Assume that $v$ is divergence-free. There exists
a constant C, depending only on $d$ and $p$, such that
$$\|\Pi(v,v)\|_{L^p}\leq C\|v\|_{L^p}\|\nabla v\|_{L^\infty}.$$
\end{lemm}
\subsection{Estimates for the electronic potential $\phi$}
In order to control $\nabla \phi,$ we introduce the following Hardy-Littlewood-Sobolev inequality.
\begin{lemm}\cite{Lemari¨¦-Rieusset}\label{Lemari¨¦-Rieusset}
For $1<p<\infty$ and $0<\gamma<\frac{d}{p},$ the operator $(-\triangle)^{\frac{\gamma}{2}}$ is bounded from $L^p$ to $L^q$ where $\frac{1}{q}=\frac{1}{p}-\frac{\gamma}{d}.$
\end{lemm}
As a consequence, we have the following lemmas.
\begin{lemm}\label{phi}
Let $s\in \mathbb{R},$ $1<p<d$ and $1\leq r\leq \infty.$ If $a\in B^{s}_{p,r},$ then \begin{align*}
  \nabla(-\triangle)^{-1}a\in L^q+B^{s+1}_{p,r},
\end{align*}
where $\frac{1}{q}=\frac{1}{p}-\frac{1}{d}.$
\end{lemm}
\noindent{Proof.} We split $\nabla(-\triangle)^{-1}a$ into low and high frequencies
\begin{align*}
\nabla(-\triangle)^{-1}a=\nabla(-\triangle)^{-1}\triangle_{-1}a
+\nabla(-\triangle)^{-1}(Id-\triangle_{-1})a.
\end{align*}
Note that $R_j=\frac{\partial_i}{(-\triangle)^{\frac{1}{2}}}$ is a Calderon-Zygmund operator. Combining with Lemma \ref{Lemari¨¦-Rieusset}, we have
\begin{align*}
\nabla(-\triangle)^{-1}\triangle_{-1}a
=\frac{\nabla}{(-\triangle)^{\frac{1}{2}}}\frac{1}
{(-\triangle)^{\frac{1}{2}}}\triangle_{-1}a
\in L^q.
\end{align*}
Next, by virtue of Lemma \ref{Bi}, we get
\begin{align*}
  \|\nabla(-\triangle)^{-1}(Id-\triangle_{-1})a\|_{B^{s+1}_{p,r}}\lesssim
  \|(Id-\triangle_{-1})a\|_{B^{s}_{p,r}}\lesssim\|a\|_{B^{s}_{p,r}}.
\end{align*}
This completes the proof of the lemma.
\begin{lemm}\label{possion}
 Under the assumption of Lemma \ref{phi},
  there exists a function $\phi\in B^{s+2}_{p,r}+L^\infty_{p^\alpha}$ for some $\alpha \in(0,1)$ depending on $p$ and d, such that
  \begin{align}
  \left\{
    \begin{array}{l}
    \triangle \phi=a, \\[1ex]
    \nabla \phi=\nabla(-\triangle)^{-1}a.
    \end{array}
    \right.
  \end{align}
\end{lemm}
\noindent{Proof.} For the high frequency part, using Lemma \ref{Bi}, we get that
\begin{align*}
\phi_0\triangleq-(-\triangle)^{-1}(Id-\triangle_{-1})a\in B^{s+2}_{p,r}
\end{align*}
solves
\begin{align*}
  \left\{
    \begin{array}{l}
    \triangle \phi_0=(Id-\triangle_{-1})a, \\[1ex]
    \nabla \phi_0=-\nabla(-\triangle)^{-1}(Id-\triangle_{-1})a.
    \end{array}
    \right.
  \end{align*}
Next for the low frequency part, by Lemma \ref{Lemari¨¦-Rieusset}, we infer that   
\begin{align*}
  f=(f_1,\cdots,f_d)\triangleq-\nabla(-\triangle)^{-1}\triangle_{-1}a\in (L^{q})^d,
\end{align*}
where $\frac{1}{q}=\frac{1}{p}-\frac{1}{d}.$\\
Again using Lemma \ref{Bi}, we have
$$f\in(W^{m,s}(\mathbb{R}^d))^d,~\forall ~m\in\mathbb{Z}^+,~q\leq s\leq \infty,$$
from which it follows that $f(y)\rightarrow 0,$ as $|y|\rightarrow\infty.$ \\
Set
\begin{align}
-\phi_{-1}
=&\sum_{k=1}^{d}\theta\partial_kE_d\ast f_k+\sum_{k=1}^{d}\int_{\mathbb{R}^n}\big((1-\theta)\partial_kE_d\big)(x-y)\theta(\frac{x-y}{\langle x\rangle})f(y)dy \\\nonumber
&+\sum_{k=1}^{d}\int_{\mathbb{R}^n}\int_0^1x_k\partial_{tx_k}[\big((1-\theta)\partial_kE_d\big)(tx-y)]
(1-\theta)(\frac{tx-y}{\langle tx\rangle})f(y)dtdy\\\nonumber
&+\sum_{k=1}^{d}\int_{\mathbb{R}^n}\int_0^1x_k\big((1-\theta)\partial_kE_d\big)(tx-y)
\partial_{tx_k}[(1-\theta)(\frac{tx-y}{\langle tx\rangle})]f(y)dtdy\\\nonumber
\triangleq&I_1+I_2+I_3+I_4,
\end{align}
with $E_d$ and $\theta$ defined before.\\
We choose $l$ to be sufficiently large such that $q\leq l<\infty, ~and~(d-1)(l'-1)<1
,$ where $l'$ is the conjugate of $l.$
 A direct calculation yields that
\begin{align*}
I_1\lesssim&\int_{|x-y|\leq 2}\frac{1}{|x-y|^{d-1}}|f(y)|dy\lesssim\|f\|_{L^\infty},
\end{align*}
\begin{align*}
I_2\lesssim&\int_{1\leq|x-y|\leq 2\langle x\rangle}\frac{1}{|x-y|^{d-1}}|f(y)|dy\lesssim\|f\|_{L^l}\langle x\rangle^{1-(d-\frac{d}{l'})},
\end{align*}
\begin{align*}
I_3\lesssim&\int_0^1\int_{|tx-y|\geq\langle tx\rangle}\frac{|x|}{|tx-y|^{d}}|f(y)|dydt
\lesssim\int_0^1\|f\|_{L^l}\frac{|x|}{{\langle tx\rangle}^{d-\frac{d}{l'}}}dt\\
\lesssim&\int_1^{2\langle x\rangle}\frac{1}{s^{d-\frac{d}{l'}}}ds\|f\|_{L^l}
\lesssim{\langle x\rangle}^{1-(d-\frac{d}{l'})}\|f\|_{L^l},
\end{align*}
\begin{align*}
I_4\lesssim&\int_0^1\int_{\langle tx\rangle\leq|tx-y|\leq2\langle tx\rangle}|x|\frac{1}{|tx-y|^{d-1}}\frac{1}{\langle tx\rangle}|f(y)|dydt
\lesssim\int_0^1 \frac{|x|}{\langle tx\rangle}{\langle tx\rangle}^{1-d+\frac{d}{l'}}\|f\|_{L^l}dt\\
\lesssim&{\langle x\rangle}^{1-(d-\frac{d}{l'})}\|f\|_{L^l},
\end{align*}
from which it follows that
$$
  \phi_{-1}\in L^\infty_{P^\alpha},
~\textit{where} ~\alpha=1-d+\frac{d}{l'}\in(0,1).$$
Note that
\begin{align*}
-\phi_{-1}=&\sum_{k=1}^{d}\theta\partial_kE_d\ast f_k+\sum_{k=1}^{d}\int_{\mathbb{R}^n}[\big((1-\theta)\partial_kE_d\big)(x-y)
-\big((1-\theta)^2\partial_kE_d\big)(-y)]f_k(y)dy.
\end{align*}
By virtue of the Lebesgue dominated convergence theorem, we have
\begin{align*}
  \nabla\phi_{-1}=-\sum_{k=1}^{d}\theta\partial_kE_d\ast \nabla f_k-\sum_{k=1}^{d}\int_{\mathbb{R}^n}[\nabla_x((1-\theta)\partial_kE_d)(x-y)
]f_k(y)dy.\end{align*}
As $\nabla\partial_k(-\triangle)^{-1}$ is a Calderon-Zygmund operator, we infer that
$\nabla\partial_kE_d\triangle_{-1}a\in L^p,$ hence
\begin{align*}
&-\sum_{k=1}^{d}\int_{\mathbb{R}^n}[\nabla_x\big((1-\theta)\partial_kE_d\big)(x-y)
]f_k(y)dy\\=&-\sum_{k=1}^{d}\int_{\mathbb{R}^n}\big((1-\theta)\partial_kE_d\big)(x-y)
\nabla_y f_k(y)dy\\
=&\sum_{k=1}^{d}\int_{\mathbb{R}^n}\big((1-\theta)\partial_kE_d\big)(x-y)
\nabla_y \partial_k(-\triangle)^{-1}\triangle_{-1}a(y)dy,
\end{align*}
which results in
\begin{align*}
\nabla\phi_{-1}=\sum_{k=1}^{d}\partial_kE_d\ast  \nabla \partial_k(-\triangle)^{-1}\triangle_{-1}a=-\nabla (-\triangle)^{-1}\triangle_{-1}a=f,\end{align*}and then \begin{align*}
\triangle\phi_{-1}=div ~\nabla\phi_{-1}=div~ f=\triangle_{-1}a.
\end{align*}
Finally, letting $\phi=\phi_{0}+\phi_{-1},$ we complete the proof of the
lemma.\qed
\section{Proof of Theorem \ref{a1}}
\hspace{0.5cm}To begin, we mention that $1<p_2<d$ and the conditions $(\ref{1.7})-(\ref{1.9})$ imply that
\begin{align}\label{jingchangyong}1<\frac{d}{p_2}< s_2+\frac{3}{2}-s_1+\frac{d}{p_1}\leq s_2+\frac{3}{2}-1\leq s_2+\frac{1}{2},\end{align} which will be frequently used.\\
We assume that $s'_1$ satisfies (\ref{s'}), $\varepsilon_1$ satisfies (\ref{ee}),
and that $q_2$ satisfies $\frac{1}{q_2}=\frac{1}{p_2}-\frac{1}{d}.$

For the existence part of Theorem \ref{a1}, we solve the $\widetilde{ENPP}$ system first, and then we show that the solution of the $\widetilde{ENPP}$ system does provide a solution for the $ENPP$ system.
  \subsection{Existence for the $\widetilde{ENPP}$ system}
  \subsubsection{First step: Construction of approximate solutions and uniform bounds}
  In order to define a sequence $(u^m,n^m,p^m)|_{m\in\mathbb{N}}$ of global approximate solutions to the $\widetilde{ENPP}$ system, we use an iterative scheme.
  First we set $u^0=u_0,~n^0=e^{t\triangle}n_0,~p^0=e^{t\triangle}p_0.$ Thanks to Lemma \ref{dts}, it is easy to see that
    $$(u^0,n^0,p^0)\in \widetilde{L}^\infty_{loc}(B^{s_1}_{p_1,r_1})\times \Big(\widetilde{L}^\infty_{loc}(B^{s_2}_{p_2,r_2})\cap\widetilde{L}^1_{loc}(B^{s_2+2}_{p_2,r_2})\Big)^2,$$
    and
   \begin{align*}
   &\|u^0\|_{\widetilde{L}^\infty_t(B^{s_1}_{p_1,r_1})}
   +\|n^0\|_{\widetilde{L}^\infty_t(B^{s_2}_{p_2,r_2})\cap\widetilde{L}^1_t(B^{s_2+2}_{p_2,r_2})}
   +\|p^0\|_{\widetilde{L}^\infty_t(B^{s_2}_{p_2,r_2})\cap\widetilde{L}^1_t(B^{s_2+2}_{p_2,r_2})}\\\leq &C(1+t)(\|u_0\|_{B^{s_1}_{p_1,r_1}}+\|n_0\|_{B^{s_2}_{p_2,r_2}}+\|p\|_{B^{s_2}_{p_2,r_2}}).
   \end{align*}
   Then, assuming that $$(u^m,n^m,p^m)\in \widetilde{L}^\infty_{loc}(B^{s_1}_{p_1,r_1})\times \Big(\widetilde{L}^\infty_{loc}(B^{s_2}_{p_2,r_2})\cap\widetilde{L}^1_{loc}({B^{s_2+2}_{p_2,r_2}})\Big)^2,$$
   we solve the following linear system:
  \begin{align}
\left\{
\begin{array}{c}
u^{m+1}_t+u^{m}\cdot \nabla u^{m+1}+\Pi(u^{m},u^{m})=\mathcal{P}\big((n-p)^{m}\psi^{m}\big),  \\[1ex]
n^{m+1}_t-\triangle n^{m+1}=-\big(T_{u^{m}}\nabla n^{m}+T_{\nabla n^{m}}u^{m}+\nabla\cdot R(u^{m}n^{m})\big)-\nabla\cdot(n^{m}\psi^{m}), \\[1ex]
p^{m+1}_t-\triangle p^{m+1}=-\big(T_{u^{m}}\nabla p^{m}+T_{\nabla p^{m}}u^{m}+\nabla\cdot R(u^{m}p^{m})\big)+\nabla\cdot(p^{m}\psi^{m}), \\[1ex]
\psi^{m}=\nabla(\triangle)^{-1}(p-n)^{m},\\[1ex]
(u^{m+1},n^{m+1},p^{m+1})|_{t=0}=(u_0,n_0,p_0).
\end{array}
\right.
\end{align}
Lemma \ref{ts} ensures that 
\begin{align}\label{u'}
  \|u^{m+1}\|_{\widetilde{L}^\infty_t(B^{s_1}_{p_1,r_1})}\lesssim & exp(C\int_0^t\|u^{m}\|_{B^{s_1}_{p_1,r_1}}dt')\Big(
  \|u_0\|_{B^{s_1}_{p_1,r_1}}\\\nonumber&+\|\Pi(u^{m},u^{m})\|_{\widetilde{L}^1_t(B^{s_1}_{p_1,r_1})}+
  \|\mathcal{P}((n^{m}-p^{m})\psi^{m})\|_{\widetilde{L}^1_t(B^{s_1}_{p_1,r_1})}\Big).
\end{align}
 Using Lemma \ref{pi}, we get
 \begin{align}\label{pii}
   \|\Pi(u,u)\|_{\widetilde{L}^1_t(B^{s_1}_{p_1,r_1})}\lesssim \|u\|_{\widetilde{L}^\infty_t(B^{s_1}_{p_1,r_1})}\|u\|_{\widetilde{L}^\infty_t(B^{s_1}_{p_1,r_1})}t,
 \end{align}
where we have used the fact that $B^{s_1}_{p_1,r_1}\hookrightarrow C^{0,1}.$\\
As for the term $\mathcal{P}\big((n-p)^{m}\psi^{m}\big),$ we first consider the case $p_2\leq p_1.$ We have
  \begin{align}\label{yi}
  &\|\mathcal{P}\big((n^{m}-p^{m})\psi^{m}\big)\|_{B^{s_1}_{p_1,r_1}}\lesssim
  \|\mathcal{P}\big((n^{m}-p^{m}\big)\psi^{m})\|_{B^{s_1-\frac{d}{p_1}+\frac{d}{p_2}}_{p_2,r_1}} \lesssim
  \|(n^m-p^m)\psi^{m}\|_{B^{s_2+\frac{3}{2}-2\varepsilon_1}_{p_2,r_1}}\\\nonumber \lesssim&
  \|T_{n^m-p^m}(Id-\triangle_{-1})\psi^{m}\|_{B^{s_2+\frac{3}{2}-2\varepsilon_1}_{p_2,r_1}}
  +\|T_{\psi^m}(n^{m}-p^{m})+R(n^m-p^m,\psi^{m})\|_{B^{s_2+\frac{3}{2}-2\varepsilon_1}_{p_2,r_1}}\\\nonumber
  \lesssim&
  \|n^{m}-p^{m}\|_{B^{-\varepsilon_1}_{\infty,r_2}}\|(Id-\triangle_{-1})\psi^{m}\|_{B^{s_2+\frac{3}{2}-\varepsilon_1}_{p_2,r_2}}
  +\|\psi^{m}\|_{B^{-\varepsilon_1}_{\infty,r_2}}\|n^{m}-p^{m}\|_{B^{s_2+\frac{3}{2}-\varepsilon_1}_{p_2,r_2}}\\ \nonumber \lesssim&\|n^{m}-p^{m}\|_{B^{\frac{d}{p_2}-\varepsilon_1}_{p_2,r_2}}\|n^{m}-p^{m}\|_{B^{s_2+\frac{1}{2}-\varepsilon_1}_{p_2,r_2}}
  +\|\psi^m\|_{B^{\frac{d}{q_2}-\varepsilon_1}_{q_2,r_2}}\|(n^{m}-p^{m})\|_{B^{s_2+\frac{3}{2}-\varepsilon_1}_{p_2,r_2}}
  \\\nonumber \lesssim&\|n^{m}-p^{m}\|_{B^{s_2+\frac{1}{2}-\varepsilon_1}_{p_2,r_2}}\|n^{m}-p^{m}\|_{B^{s_2+\frac{1}{2}-\varepsilon_1}_{p_2,r_2}}
  +\|n^{m}-p^{m}\|_{B^{\frac{d}{p_2}-1-\varepsilon_1}_{p_2,r_2}}\|n^{m}-p^{m}\|_{B^{s_2+\frac{3}{2}-\varepsilon_1}_{p_2,r_2}}
  \\\nonumber\lesssim&\|n^{m}-p^{m}\|_{B^{s_2+\frac{1}{2}-\varepsilon_1}_{p_2,r_2}}\|n^{m}-p^{m}\|_{B^{s_2+\frac{1}{2}-\varepsilon_1}_{p_2,r_2}}
  +\|n^{m}-p^{m}\|_{B^{s_2-\frac{1}{2}-\varepsilon_1}_{p_2,r_2}}\|n^{m}-p^{m}\|_{B^{s_2+\frac{3}{2}-\varepsilon_1}_{p_2,r_2}}
  \\\nonumber \lesssim&
  \|n^{m}-p^{m}\|_{B^{s_2+\frac{1}{2}-\varepsilon_1}_{p_2,r_2}}\|n^{m}-p^{m}\|_{B^{s_2+\frac{3}{2}-\varepsilon_1}_{p_2,r_2}}.
  \end{align}
  In the case $p_2>p_1,$ let $\frac{1}{p_3}=\frac{1}{p_1}-\frac{1}{p_2}.$ The condition $(\ref{1.9})$ implies that $\frac{1}{p_3}\leq\frac{1}{q_2}\leq\frac{1}{p_2}.$ Then we have
  \begin{align*}
  &\|\mathcal{P}\big((n^{m}-p^{m})\psi^{m}\big)\|_{B^{s_1}_{p_1,r_1}}\lesssim
  \|(n^{m}-p^{m})\psi^{m}\|_{B^{s_1}_{p_1,r_1}}\\\lesssim&
  \|n^{m}-p^{m}\|_{B^{-\varepsilon_1}_{p_3,r_2}}\|(Id-\triangle_{-1})\psi^{m}\|_{B^{s_1+\varepsilon_1}_{p_2,r_2}}
  +\|\psi^{m}\|_{B^{-\varepsilon_1}_{p_3,r_2}}\|n^{m}-p^{m}\|_{B^{s_1+\varepsilon_1}_{p_2,r_2}}\\ \lesssim&\|n^{m}-p^{m}\|_{B^{\frac{d}{p_2}-\frac{d}{p_3}-\varepsilon_1}_{p_2,r_2}}\|n^{m}-p^{m}\|_{B^{s_1-1+\varepsilon_1}_{p_2,r_2}}
  +\|\psi^{m}\|_{B^{\frac{d}{q_2}-\frac{d}{p_3}-\varepsilon_1}_{q_2,r_2}}\|n^{m}-p^{m}\|_{B^{s_2+\frac{3}{2}-\varepsilon_1}_{p_2,r_2}}
  \\\lesssim&\|n^{m}-p^{m}\|_{B^{s_2+\frac{1}{2}-\varepsilon_1}_{p_2,r_2}}\|n^{m}-p^{m}\|_{B^{s_2+\frac{1}{2}-\varepsilon_1}_{p_2,r_2}}
  +\|n^{m}-p^{m}\|_{B^{\frac{d}{p_2}-1-\frac{d}{p_3}-\varepsilon_1}_{p_2,r_2}}\|n^m-p^{m}\|_{B^{s_2+\frac{3}{2}-\varepsilon_1}_{p_2,r_2}}
  \\\lesssim&\|n^{m}-p^{m}\|_{B^{s_2+\frac{1}{2}-\varepsilon_1}_{p_2,r_2}}\|n^{m}-p^{m}\|_{B^{s_2+\frac{1}{2}-\varepsilon_1}_{p_2,r_2}}
  +\|n^{m}-p^{m}\|_{B^{s_2-\frac{1}{2}-\varepsilon_1}_{p_2,r_2}}\|n^{m}-p^{m}\|_{B^{s_2+\frac{3}{2}-\varepsilon_1}_{p_2,r_2}}
  \\\lesssim&
  \|n^{m}-p^{m}\|_{B^{s_2+\frac{1}{2}-\varepsilon_1}_{p_2,r_2}}\|n^{m}-p^{m}\|_{B^{s_2+\frac{3}{2}-\varepsilon_1}_{p_2,r_2}}.
  \end{align*}
  Hence, combining the above two estimates, we obtain
  \begin{align}\label{np}
    \|\mathcal{P}\big(n^{m}-p^{m}\psi^{m}\big)\|_{\widetilde{L}^1(B^{s_1}_{p_1,r_1})}\lesssim
    \|n^{m}-p^{m}\|_{\widetilde{L}^{\frac{4}{1-2\varepsilon_1}}(B^{s_2+\frac{1}{2}-\varepsilon_1}_{p_2,r_2})}
    \|n^{m}-p^{m}\|_{\widetilde{L}^{\frac{4}{3-2\varepsilon_1}}(B^{s_2+\frac{3}{2}-\varepsilon_1}_{p_2,r_2})}t^{\varepsilon_1}.
  \end{align}
  Inserting this inequality and (\ref{pii}) into (\ref{u'}), we get
\begin{align}\label{u}
  \|u^{m+1}\|_{\widetilde{L}^\infty_t(B^{s_1}_{p_1,r_1})}\lesssim & exp(C\int_0^t\|u^{m}\|_{B^{s_1}_{p_1,r_1}}dt')\Big(
  \|u_0\|_{B^{s_1}_{p_1,r_1}}+\|u\|_{\widetilde{L}^\infty_t(B^{s_1}_{p_1,r_1})}^2t
  \\\nonumber&+\|n^{m}-p^{m}\|_{\widetilde{L}^{\frac{4}{3-2\varepsilon_1}}(B^{s_2+\frac{1}{2}-\varepsilon_1}_{p_2,r_2})}
    \|n^{m}-p^{m}\|_{\widetilde{L}^{\frac{4}{1-2\varepsilon_1}}(B^{s_2+\frac{3}{2}-\varepsilon_1}_{p_2,r_2})}t^{\varepsilon_1}\Big).
\end{align}

As regards $n^{m+1},$ it follows from Lemma \ref{dts} that
\begin{align*}
  &\|n^{m+1}\|_{\widetilde{L}^\infty_t(B^{s_2}_{p_2,r_2})}
  +\|n^{m+1}\|_{\widetilde{L}^1_t(B^{s_2+2}_{p_2,r_2})}\\\lesssim
  &(1+t)\Big(\|n_0\|_{B^{s_2}_{p_2,r_2}}+\|T_{u^{m}}\nabla n^{m}\|_{\widetilde{L}^1_t(B^{s_2}_{p_2,r_2})}
  +\|T_{\nabla n^{m}}u^{m}\|_{\widetilde{L}^1_t(B^{s_2}_{p_2,r_2})}\\&+\|\nabla\cdot R(u^{m}n^{m})\|_{\widetilde{L}^1_t(B^{s_2}_{p_2,r_2})}
  +\|\nabla\cdot(n^{m}\psi^m)\|_{\widetilde{L}^1_t(B^{s_2}_{p_2,r_2})}\Big).
  \end{align*}
  According to Lemmas \ref{T} and \ref{R}, we get
  \begin{align}\label{1}
   &\|T_{u^{m}}\nabla n^{m}\|_{\widetilde{L}^1_t(B^{s_2}_{p_2,r_2})}+
   \|\nabla\cdot R(u^{m}n^{m}) \|_{\widetilde{L}^1_t(B^{s_2}_{p_2,r_2})}
   \\\nonumber\lesssim&
   \|u^{m}\|_{L^\infty_t(L^\infty)}\|n^{m}\|_{\widetilde{L}^1_t(B^{s_2+1}_{p_2,r_2})}
   \lesssim\|u^{m}\|_{\widetilde{L}^\infty_t(B^{s_1}_{p_1,r_1})}
   \|n^{m}\|_{\widetilde{L}^2_t(B^{s_2+1}_{p_2,r_2})}t^{\frac{1}{2}},
 \end{align}
 \begin{align}\label{2}
   &\|\nabla (n^m\psi^m)\|_{\widetilde{L}^1_t(B^{s_2}_{p_2,r_2})}\lesssim\| (n^m\psi^m)\|_{\widetilde{L}^1_t(B^{s_2+1}_{p_2,r_2})}
  \\\nonumber\lesssim&
\|n^m\|_{L^1_t(L^\infty)}
   \|(Id-\triangle_{-1})\psi^m\|_{\widetilde{L}^\infty_t(B^{s_2+1}_{p_2,r_2})}+ \|\psi^m\|_{\widetilde{L}^\infty_t(B^{\frac{d}{q_2}}_{q_2,1})}
   \|n^m\|_{\widetilde{L}^1_t(B^{s_2+1}_{p_2,r_2})}\\
   \nonumber\lesssim&\|n^m\|_{\widetilde{L}^1_t(B^{s_2+\frac{1}{2}}_{p_2,r_2})}
   \|n^m-p^m\|_{\widetilde{L}^\infty_t(B^{s_2}_{p_2,r_2})}+
   \|n^m-p^m\|_{\widetilde{L}^\infty_t(B^{\frac{d}{p_2}-1}_{p_2,1})}
  \|n^m\|_{L^1_t(B^{s_2+1}_{p_2,r_2})}\\
  \nonumber\lesssim&\|n^m\|_{\widetilde{L}^1_t(B^{s_2+\frac{1}{2}}_{p_2,r_2})}
   \|n^m-p^m\|_{\widetilde{L}^\infty_t(B^{s_2}_{p_2,r_2})}+
   \|n^m-p^m\|_{\widetilde{L}^\infty_t(B^{s_2-\frac{1}{2}}_{p_2,r_2})}
   \|n^m\|_{\widetilde{L}^1_t(B^{s_2+1}_{p_2,r_2})}\\
   \nonumber\lesssim&\|n^m-p^m\|_{\widetilde{L}^\infty_t(B^{s_2}_{p_2,r_2})}
   \|n^m\|_{\widetilde{L}^2_t(B^{s_2+1}_{p_2,r_2})}t^{\frac{1}{2}}.
 \end{align}
 In order to bound $T_{\nabla n^{m}}u^{m},$ we first consider the case $p_2>p_1.$ As $s_2+\frac{d}{p_1}-\frac{d}{p_2}-s_1<0,$ we have \begin{align*}
   \|T_{\nabla n^{m}}u^{m}\|_{B^{s_2}_{p_2,r_2}}\lesssim&
   \|\nabla n^m\|_{B^{s_2+\frac{d}{p_1}-\frac{d}{p_2}-s_1}_{\infty,r_2}}
   \|u^m\|_{B^{s_1-\frac{d}{p_1}+\frac{d}{p_2}}_{p_2,r_1}}
   \\ \lesssim&\|\nabla n^m\|_{B^{s_2+\frac{d}{p_1}-s_1}_{p_2,r_2}}
   \|u^m\|_{B^{s_1}_{p_1,r_1}}\lesssim\| n^m\|_{B^{s_2}_{p_2,r_2}}
   \|u^m\|_{B^{s_1}_{p_1,r_1}}.
 \end{align*}
 The case $p_2\leq p_1$ works in almost the same way. Letting $\frac{1}{p_3}=\frac{1}{p_2}-\frac{1}{p_1},$ we have
 \begin{align*}
   \|T_{\nabla n^{m}}u^{m}\|_{B^{s_2}_{p_2,r_2}}\lesssim&
   \|\nabla n^m\|_{B^{s_2-s_1}_{p_3,r_2}}
   \|u^m\|_{B^{s_1}_{p_1,r_1}}
   \\ \lesssim&\|\nabla n^m\|_{B^{s_2-\frac{d}{p_3}+\frac{d}{p_2}-s_1}_{p_2,r_2}}
   \|u^m\|_{B^{s_1}_{p_1,r_1}}\lesssim\| n^m\|_{B^{s_2}_{p_2,r_2}}
   \|u^m\|_{B^{s_1}_{p_1,r_1}}.
 \end{align*}
Hence, we obtain
\begin{align}\label{3}
  \|T_{\nabla n^{m}}u^{m}\|_{\widetilde{L}^1_t(B^{s_2}_{p_2,r_2})}\lesssim\| n^m\|_{\widetilde{L}^\infty_t(B^{s_2}_{p_2,r_2})}
   \|u^m\|_{\widetilde{L}^\infty_t(B^{s_1}_{p_1,r_1})}t.
\end{align}
 Thus, we conclude that
 \begin{align}\label{n}
  &\|n^{m+1}\|_{\widetilde{L}^\infty_t(B^{s_2}_{p_2,r_2})}
  +\|n^{m+1}\|_{\widetilde{L}^1_t(B^{s_2+2}_{p_2,r_2})}\\\nonumber\lesssim
  &(1+t)\Big(\|n_0\|_{B^{s_2}_{p_2,r_2}}+\|u^{m}\|_{\widetilde{L}^\infty_t(B^{s_1}_{p_1,r_1})}
   \|n^{m}\|_{\widetilde{L}^2_t(B^{s_2+1}_{p_2,r_2})}t^{\frac{1}{2}}
  \\\nonumber&+\|n^m-p^m\|_{\widetilde{L}^\infty_t(B^{s_2}_{p_2,r_2})}
   \|n^m\|_{\widetilde{L}^2_t(B^{s_2+1}_{p_2,r_2})}t^{\frac{1}{2}}
  +\|n^m\|_{\widetilde{L}^\infty_t(B^{s_2}_{p_2,r_2})}
   \|u^m\|_{\widetilde{L}^\infty_t(B^{s_1}_{p_1,r_1})}t\Big).
  \end{align}

  A similar process as above ensures that
  \begin{align}\label{p}
  &\|p^{m+1}\|_{\widetilde{L}^\infty_t(B^{s_2}_{p_2,r_2})}
  +\|p^{m+1}\|_{\widetilde{L}^1_t(B^{s_2+2}_{p_2,r_2})}\\\nonumber\lesssim
  &(1+t)\Big(\|p_0\|_{B^{s_2}_{p_2,r_2}}+\|u^{m}\|_{\widetilde{L}^\infty_t(B^{s_1}_{p_1,r_1})}
   \|p^{m}\|_{\widetilde{L}^2_t(B^{s_2+1}_{p_2,r_2})}t^{\frac{1}{2}}
  \\\nonumber&+\|n^m-p^m\|_{\widetilde{L}^\infty_t(B^{s_2}_{p_2,r_2})}
   \|p^m\|_{\widetilde{L}^2_t(B^{s_2+1}_{p_2,r_2})}t^{\frac{1}{2}}
  +\|p^m\|_{\widetilde{L}^\infty_t(B^{s_2}_{p_2,r_2})}
   \|u^m\|_{\widetilde{L}^\infty_t(B^{s_1}_{p_1,r_1})}t\Big).
  \end{align}

  Denote $$E^{m}(t)\triangleq \|u^{m}\|_{\widetilde{L}^\infty_t(B^{s_1}_{p_1,r_1})}+\|n^{m}\|_{\widetilde{L}^\infty_t(B^{s_2}_{p_2,r_2})\cap \widetilde{L}^1_t(B^{s_2+2}_{p_2,r_2})}+
  \|p^{m}\|_{\widetilde{L}^\infty_t(B^{s_2}_{p_2,r_2})\cap \widetilde{L}^1_t(B^{s_2+2}_{p_2,r_2})},$$
  and
  $$E^0\triangleq\|u_0\|_{B^{s_1}_{p_1,r_1}}+\|n_0\|_{B^{s_2}_{p_2,r_2}}+
  \|p_0\|_{B^{s_2}_{p_2,r_2}}.$$
  By using interpolation and plugging the inequalities (\ref{n}) and (\ref{p}) into (\ref{u}) yield
  \begin{align*}
    E^{m+1}(t)\leq C\big(e^{CE^{m}(t)t}+1+t\big)\Big(E^0+\big(E^{m}(t)\big)^2\big(t+t^{\frac{1}{2}}+t^{\varepsilon_1}\big)\Big).
  \end{align*}
  Let us choose a positive $T_0\leq1$ such that $exp{(8C^2E^0T_0)}\leq 2$ and $T_0^{\varepsilon_1}\leq \frac{1}{192C^2E_0}.$ The induction hypothesis then implies that $$E^{m}(T_0)\leq 8CE^0.$$
\subsubsection{Second step: Convergence of the sequence}
  Let us fix some positive $T$ such that $T\leq T_0,$ and $(2CE^0)^4T\leq 1.$ We assume that $s'$ satisfies ($\ref{s'}$).

  By taking the difference between the equations for $u^{m+1}$ and $u^{m},$ one finds that
   \begin{align}\label{uu}
     &(u^{m+1}-u^m)_t+u^{m}\cdot \nabla (u^{m+1}-u^m)\\=\nonumber&(u^{m-1}-u^{m})\nabla u^{m}-\Pi(u^{m}-u^{m-1},u^{m}+u^{m-1})\\\nonumber&+\mathcal{P}\big((n^m-p^{m})(\psi^{m}-\psi^{m-1})\big)
   +\mathcal{P}\big((n^{m}-p^{m}-n^{m-1}+p^{m-1})\psi^{m-1}\big).
   \end{align}
   Thanks to Lemmas (\ref{T}), (\ref{R}), (\ref{pi}) and (\ref{deng}), for a.e. $t\in[0,T],$ we have
   \begin{align}\label{i}
     &\|(u^{m-1}-u^{m})\nabla u^{m+1}\|_{\widetilde{L}^1_t(B^{s'_1}_{p_1,r_1})}+\|\Pi(u^{m}-u^{m-1},u^{m}+u^{m-1})\|_{\widetilde{L}^1_t(B^{s'_1}_{p_1,r_1})}
     \\\nonumber\lesssim &\|u^{m-1}-u^{m}\|_{\widetilde{L}^\infty_t(B^{s'_1}_{p_1,r_1})}
     (\|u^{m}\|_{\widetilde{L}^\infty_t(B^{s_1}_{p_1,r_1})}
     +\|u^{m-1}\|_{\widetilde{L}^\infty_t(B^{s_1}_{p_1,r_1})})t.
   \end{align}
   The process of dealing with the terms $\mathcal{P}\big((n^m-p^{m})(\psi^{m}-\psi^{m-1})\big)
   ~and~\mathcal{P}\big((n^{m}-p^{m}-n^{m-1}+p^{m-1})\psi^{m-1}\big)$ is similar as that of the inequality (\ref{np}), so we get
   \begin{align}\label{ii}
    &\|\mathcal{P}\big((n^m-p^{m})(\psi^{m}-\psi^{m-1})\big)\|_{\widetilde{L}^\infty_t(B^{s'_1}_{p_1,r_1})}\\\nonumber\lesssim
    &\|n^m-p^{m}\|_{\widetilde{L}^{4}_t(B^{s_2+\frac{1}{2}}_{p_2,r_2})}
    \|n^{m}-p^{m}-n^{m-1}+p^{m-1}\|_{\widetilde{L}^{4}_t(B^{s_2-\frac{1}{2}}_{p_2,r_2})}t^{\frac{1}{2}},
    \end{align}
     \begin{align}\label{iii}
     &\|\mathcal{P}\big((n^{m}-p^{m}-n^{m-1}+p^{m-1})\psi^{m-1}\big)\|_{\widetilde{L}^\infty_t(B^{s'_1}_{p_1,r_1})}\\
     \nonumber&\lesssim
    \|n^{m}-p^{m}-n^{m-1}+p^{m-1}\|_{\widetilde{L}^{\frac{4}{3}}_t(B^{s_2+\frac{1}{2}}_{p_2,r_2})}
    \|n^{m-1}-p^{m-1}\|_{\widetilde{L}^{\infty}_t(B^{s_2-\frac{1}{2}}_{p_2,r_2})}t^{\frac{1}{4}}\\
    \nonumber&\lesssim
    \|n^{m}-p^{m}-n^{m-1}+p^{m-1}\|_{\widetilde{L}^{\frac{4}{3}}_t(B^{s_2+\frac{1}{2}}_{p_2,r_2})}
    \|n^{m-1}-p^{m-1}\|_{\widetilde{L}^{\infty}_t(B^{s_2}_{p_2,r_2})}t^{\frac{1}{4}}.
   \end{align}
   Applying Lemma \ref{ts} to (\ref{uu}), we thus obtain
   \begin{align}\label{uug}
  &\|u^{m+1}-u^{m}\|_{\widetilde{L}^\infty_t(B^{s'_1}_{p_1,r_1})}\\\nonumber\lesssim & exp(C\int_0^t\|u^{m}\|_{B^{s_1}_{p_1,r_1}}dt')\Big(\|u^{m-1}-u^{m}\|_{\widetilde{L}^\infty_t(B^{s'_1}_{p_1,r_1})}
  (\|u^{m}\|_{\widetilde{L}^\infty_t(B^{s_1}_{p_1,r_1})}
     +\|u^{m-1}\|_{\widetilde{L}^\infty_t(B^{s_1}_{p_1,r_1})})t\\\nonumber&~~~+
  \|n^m-p^m\|_{\widetilde{L}^4_t(B^{s_2+\frac{1}{2}}_{p_2,r_2})}
    \|n^{m}-p^{m}-n^{m-1}+p^{m-1}\|_{\widetilde{L}^4_t(B^{s_2-\frac{1}{2}}_{p_2,r_2})}t^{\frac{1}{2}}
   \\\nonumber& ~~~+\|n^{m}-p^{m}-n^{m-1}+p^{m-1}\|_{\widetilde{L}^{\frac{4}{3}}_t(B^{s_2+\frac{1}{2}}_{p_2,r_2})}
    \|n^{m-1}-p^{m-1}\|_{\widetilde{L}^\infty_t(B^{s_2}_{p_2,r_2})}t^{\frac{1}{4}}\Big).
\end{align}

Note that
\begin{align*}
  (n^{m+1}-n^{m})_t-&\triangle \Big(n^{m+1}-n^{m})=-(T_{(u^{m}-u^{m-1})}\nabla n^{m}+T_{u^{m-1}}\nabla (n^{m}-n^{m-1})\\&+T_{\nabla (n^{m}-n^{m-1})}u^{m-1}+T_{\nabla n^{m}}(u^{m}-u^{m-1})+\nabla\cdot R\big((u^{m}-u^{m-1}),n^{m}\big)\\&+\nabla\cdot R(u^{m-1},n^{m}-n^{m-1})\Big)-\nabla\cdot\big((n^{m}-n^{m-1})\big)\psi^{m}+n^{m-1}(\psi^{m}-\psi^{m-1}).
\end{align*}
Following along almost the same lines of the proof of the inequalities (\ref{1})-(\ref{3}), we get
\begin{align}\label{44}
   &\|T_{(u^{m}-u^{m-1})}\nabla n^{m}+T_{\nabla n^{m}}(u^{m}-u^{m-1})+\nabla \cdot R\big((u^{m}-u^{m-1}),n^{m}\big)
    \|_{\widetilde{L}^1_t(B^{s_2-1}_{p_2,r_2})}
   \\\nonumber\lesssim&
   \|u^{m}-u^{m-1}\|_{\widetilde{L}^\infty_t(B^{s'_1}_{p_1,r_1})}
   \|n^{m}\|_{\widetilde{L}^\infty_t(B^{s_2}_{p_2,r_2})}t,
   \end{align}
   \begin{align}\label{55}
   &\|T_{u^{m}}\nabla (n^{m}-n^{m-1})+T_{\nabla (n^{m}-n^{m-1})}u^{m-1}+\nabla \cdot R(u^{m-1},(n^{m}-n^{m-1}))
    \|_{\widetilde{L}^1_t(B^{s_2-1}_{p_2,r_2})}
   \\\nonumber\lesssim&
   \|u^{m}\|_{\widetilde{L}^\infty_t(B^{s'_1}_{p_1,r_1})}
   \|n^{m}-n^{m-1}\|_{\widetilde{L}^1_t(B^{s_2}_{p_2,r_2})}\lesssim
   \|u^{m}\|_{\widetilde{L}^\infty_t(B^{s_1}_{p_1,r_1})}
   \|n^{m}-n^{m-1}\|_{\widetilde{L}^2_t(B^{s_2}_{p_2,r_2})}t^{\frac{1}{2}},
   \end{align}
   \begin{align}\label{66}
   \|\nabla\cdot\big((n^{m}-n^{m-1})\psi^{m}\big)
    \|_{\widetilde{L}^1_t(B^{s_2-1}_{p_2,r_2})}
    \lesssim&
   \|n^{m}-n^{m-1}\|_{\widetilde{L}^1_t(B^{s_2}_{p_2,r_2})}
   \|n^{m}-p^m\|_{\widetilde{L}^\infty_t(B^{s_2-\frac{1}{2}}_{p_2,r_2})}\\\nonumber\lesssim&
   \|n^{m}-n^{m-1}\|_{\widetilde{L}^2_t(B^{s_2}_{p_2,r_2})}
   \|n^{m}-p^m\|_{\widetilde{L}^\infty_t(B^{s_2}_{p_2,r_2})}t^{\frac{1}{2}},
   \end{align}
   \begin{align}\label{77}
   \|\nabla\cdot\big(n^{m-1}(\psi^{m}-\psi^{m-1})\big)
    \|_{\widetilde{L}^1_t(B^{s_2-1}_{p_2,r_2})}
  \lesssim&
   \|n^{m-1}\|_{\widetilde{L}^\infty_t(B^{s_2}_{p_2,r_2})}
   \|n^{m}-p^{m}-n^{m-1}+p^{m-1}\|_{\widetilde{L}^1_t(B^{s_2-\frac{1}{2}}_{p_2,r_2})}
   \\\nonumber\lesssim&
   \|n^{m-1}\|_{\widetilde{L}^\infty_t(B^{s_2}_{p_2,r_2})}
   \|n^{m}-p^{m}-n^{m-1}+p^{m-1}\|_{\widetilde{L}^{4}_t(B^{s_2-\frac{1}{2}}_{p_2,r_2})}t^{\frac{3}{4}}.
 \end{align}
 Hence Lemma \ref{dts} implies that
 \begin{align}\label{nng}
  &\|n^{m+1}-n^{m}\|_{\widetilde{L}^\infty_t(B^{s_2-1}_{p_2,r_2})}
  +\|n^{m+1}-n^{m}\|_{\widetilde{L}^1_t(B^{s_2+1}_{p_2,r_2})}\\\nonumber\lesssim
  &(1+t)\Big(\|u^{m}-u^{m-1}\|_{\widetilde{L}^\infty_t(B^{s'_1}_{p_1,r_1})}
   \|n^{m}\|_{\widetilde{L}^\infty_t(B^{s_2}_{p_2,r_2})}t+
  \|u^{m}\|_{\widetilde{L}^\infty_t(B^{s_1}_{p_1,r_1})}
   \|n^{m}-n^{m-1}\|_{\widetilde{L}^2_t(B^{s_2}_{p_2,r_2})}t^{\frac{1}{2}}
  \\\nonumber&+\|n^{m}-n^{m-1}\|_{\widetilde{L}^2_t(B^{s_2}_{p_2,r_2})}
   \|n^{m}-p^m\|_{\widetilde{L}^\infty_t(B^{s_2}_{p_2,r_2})}t^{\frac{1}{2}}\\\nonumber&+\|n^{m-1}\|_{\widetilde{L}^\infty_t(B^{s_2}_{p_2,r_2})}
   \|n^{m}-p^{m}-n^{m-1}+p^{m-1}\|_{\widetilde{L}^{4}_t(B^{s_2-\frac{1}{2}}_{p_2,r_2})}t^{\frac{3}{4}}\Big).
  \end{align}
  Similarly, we get
  \begin{align}\label{ppg}
  &\|p^{m+1}-p^{m}\|_{\widetilde{L}^\infty_t(B^{s_2-1}_{p_2,r_2})}
  +\|p^{m+1}-p^{m}\|_{\widetilde{L}^1_t(B^{s_2+1}_{p_2,r_2})}\\\nonumber\lesssim
  &(1+t)\Big(\|u^{m}-u^{m-1}\|_{\widetilde{L}^\infty_t(B^{s'_1}_{p_1,r_1})}
   \|p^{m}\|_{\widetilde{L}^\infty_t(B^{s_2}_{p_2,r_2})}t+
  \|u^{m}\|_{\widetilde{L}^\infty_t(B^{s_1}_{p_1,r_1})}
   \|p^{m}-p^{m-1}\|_{\widetilde{L}^2_t(B^{s_2}_{p_2,r_2})}t^{\frac{1}{2}}
  \\\nonumber&+\|p^{m}-p^{m-1}\|_{\widetilde{L}^2_t(B^{s_2}_{p_2,r_2})}
   \|n^{m}-p^m\|_{\widetilde{L}^\infty_t(B^{s_2}_{p_2,r_2})}t^{\frac{1}{2}}\\\nonumber&+\|p^{m-1}\|_{\widetilde{L}^\infty_t(B^{s_2}_{p_2,r_2})}
   \|n^{m}-p^{m}-n^{m-1}+p^{m-1}\|_{\widetilde{L}^{4}_t(B^{s_2-\frac{1}{2}}_{p_2,r_2})}t^{\frac{3}{4}}\Big).
  \end{align}
   Denote \begin{align*}F^{m}(t)\triangleq \|u^{m+1}-u^{m}\|_{\widetilde{L}^\infty_t(B^{s'_1}_{p_2,r_2})}&+\|n^{m+1}-n^{m}\|_{\widetilde{L}^\infty_t
   (B^{s_2-1}_{p_2,r_2})\cap\widetilde{L}^1_t(B^{s_2+1}_{p_2,r_2})}\\&+
  \|p^{m+1}-p^{m}\|_{\widetilde{L}^\infty_t(B^{s_2-1}_{p_2,r_2})
  \cap\widetilde{L}^1_t(B^{s_2+1}_{p_2,r_2})}.\end{align*}
  Plugging the inequalities (\ref{nng}) and (\ref{ppg}) into (\ref{uug}) yields
  \begin{align*}
    F^{m+1}(T)\leq & C\big(e^{CE^{m}(T)T}+1+T\big)\big(E^{m}(T)+E^{m-1}(T)\big)F^{m}(T)\big(T+T^{\frac{1}{2}}+T^{\frac{3}{4}}+T^{\frac{1}{4}}\big)
    \\\leq &CE^0T^{\frac{1}{4}}F^{m}(T)\leq \frac{1}{2}F^{m}(T).
  \end{align*}
 Hence, $(u^m,n^m,p^m)|_{m\in \mathbb{N}}$ is a Cauchy sequence in  $\widetilde{L}^\infty_T(B^{s'_1}_{p_1,r_1})\times \Big(\widetilde{L}^\infty_T(B^{s_2-1}_{p_2,r_2})\cap\widetilde{L}^1_T(B^{s_2+1}_{p_2,r_2})\Big)^2.$
  \subsubsection{Third step: Passing to the limit}
  Let $(u,n,p)$ be the limit of the sequence $(u^m,n^m,p^m)|_{m\in \mathbb{N}}.$ We see that $(u,n,p)\in\widetilde{L}^\infty_T(B^{s'_1}_{p_1,r_1})\times \Big(\widetilde{L}^\infty_T(B^{s_2-1}_{p_2,r_2})\cap\widetilde{L}^1_T(B^{s_2+1}_{p_2,r_2})\Big)^2.$ Using Lemma \ref{Fadou} with the uniform bounds given in Step 1, we see that $(u,n,p)\in\widetilde{L}^\infty_T(B^{s_1}_{p_1,r_1})\times \Big(\widetilde{L}^\infty_T(B^{s_2}_{p_2,r_2})\Big)^2.$ Next, by interpolating we discover that $(u^m,n^m,p^m)$ tends to $(u,n,p)$ in every space $\widetilde{L}^\infty_T(B^{s_1-\varepsilon}_{p_1,r_1})\times \Big(\widetilde{L}^\infty_T(B^{s_2-\varepsilon}_{p_2,r_2})\cap\widetilde{L}^1_T(B^{s_2+1}_{p_2,r_2})\Big)^2,$ with $\varepsilon>0,$ which suffices to pass to the limit in the $\widetilde{ENPP}$ system.

  We still have to prove that $(n,p)\in \big(\widetilde{L}^1_T(B^{s_2+2}_{p_2,r_2})\big)^2.$ In fact, it is easy to check  that $\partial_t n-\triangle n\in\widetilde{L}^1_T(B^{s_2}_{p_2,r_2}).$ Hence according to Lemma \ref{dts}, $n\in \widetilde{L}^1_T(B^{s_2+2}_{p_2,r_2}).$ Similarly, $p\in \widetilde{L}^1_T(B^{s_2+2}_{p_2,r_2}).$

 \subsection{Existence for the $ENPP$ system}
  Suppose that $(u,n,p)$ satisfies the $\widetilde{ENPP}$ system in $X(t).$  We first check that $u$ is divergence free. This may be achieved by applying $div$ to the first equation of $\widetilde{ENPP}$. We get
  \begin{align*}
    (\partial_t+u\cdot \nabla)div~u=-div~\Pi(u,u)-tr(Dv)^2.
  \end{align*}
  Lemma \ref{ts} and Lemma \ref{yyy} ensure that
  \begin{align*}
    \|div~u\|_{B^{s_1'}_{p_1,r_1}}\lesssim&
    \int_0^texp\Big(C\int_{t'}^t\|u\|_{B^{s}_{p_1,r_1}}dt''\Big)\|div~\Pi(u,u)+tr(Dv)^2\|_{B^{s_1}_{p_1,r_1}}dt'\\ \lesssim&
    \int_0^texp\Big(C\int_{t'}^t\|u\|_{B^{s_1}_{p_1,r_1}}dt''\Big)\|div~u\|_{B^{0}_{\infty,\infty}}\|u\|_{B^{s'_1+1}_{p_1,r_1}}dt',\\
    \lesssim&
    \int_0^texp\Big(C\int_{t'}^t\|u\|_{B^{s_1}_{p_1,r_1}}dt''\Big)\|div~u\|_{B^{s_1'}_{p_1,r_1}}\|u\|_{B^{s_1}_{p_1,r_1}}dt',
  \end{align*}
  where we have used $B^{s_1'}_{p_1,r_1}\hookrightarrow B^{0}_{\infty,\infty},$ and $B^{s_1'+1}_{p_1,r_1} \hookrightarrow B^{s_1}_{p_1,r_1}.$
  Using Gronwall's inequality, we conclude that $div ~u=0.$

 Next according to Lemma \ref{possion}, there exists a function $\phi\in L_T^\infty( B^{s+2}_{p,r}+L^\infty_{p^\alpha})$ for some $\alpha\in(0,1)$ satisfying
  \begin{align}\label{possion1}
  \left\{
    \begin{array}{l}
    \triangle \phi=n-p, \\[1ex]
    \nabla \phi=\nabla(-\triangle)^{-1}(p-n).
    \end{array}
    \right.
  \end{align}
  As the condition $(\ref{1.8})$ implies $s_2-\frac{d}{p_2}+2>s-\frac{d}{p_1}+\frac{1}{2}>0,$ we have
  $$\phi\in L_T^\infty(B^{s_2+2}_{p_2,r_2}+L^\infty_{p^\alpha})\hookrightarrow
  L^\infty_T(L^\infty+L^\infty_{p^\alpha})\hookrightarrow
  L^\infty_T(L^\infty_{p^\alpha}).$$

  Let $$P=P_{\pi}(u,u)-(-\triangle)^{-1}div(\triangle\phi\nabla\phi),$$
  where $P_{\pi}(u,u)\in L_T^\infty(B^{s_1+1}_{p_1,r_1}+L^\infty_L)$ is defined as in Lemma \ref{pi}.\\
  A similar argument as that of (\ref{yi}) implies that
  \begin{align*}
    \|\triangle\phi\nabla\phi\|_{B^{s_2+\frac{3}{2}-2\varepsilon_1}_{p_2,r_2}}\lesssim
    \|n-p\|_{B^{s_2+\frac{1}{2}-\varepsilon_1}_{p_2,r_2}}\|n-p\|_{B^{s_2+\frac{3}{2}-\varepsilon_1}_{p_2,r_2}}.
  \end{align*}
  Hence, by virtue of the Minkowski inequality and the imbedding inequality, we have
  \begin{align*}
    \|\triangle\phi\nabla\phi\|_{L^1_T(B^{s_2+\frac{3}{2}-2\varepsilon_1}_{p_2,r_2})}\lesssim&
    \|n-p\|_{L^{\frac{4}{3}}_T(B^{s_2+\frac{1}{2}-\varepsilon_1}_{p_2,r_2})}\|n-p\|_{L^{4}_T(B^{s_2+\frac{3}{2}-\varepsilon_1}_{p_2,r_2})}\\
    \lesssim&\|n-p\|_{L^{\frac{4}{3}}_T(B^{s_2+\frac{1}{2}-\varepsilon_1}_{p_2,1})}
    \|n-p\|_{L^{4}_T(B^{s_2+\frac{3}{2}-\varepsilon_1}_{p_2,1})}\\
  \lesssim&\|n-p\|_{\tilde{L}^{\frac{4}{3}}_T(B^{s_2+\frac{1}{2}-\varepsilon_1}_{p_2,1})}
    \|n-p\|_{\tilde{L}^{4}_T(B^{s_2+\frac{3}{2}-\varepsilon_1}_{p_2,1})}\\
    \lesssim&\|n-p\|_{\tilde{L}^{\frac{4}{3}}_T(B^{s_2+\frac{1}{2}}_{p_2,r_2})}
    \|n-p\|_{\tilde{L}^{4}_T(B^{s_2+\frac{3}{2}}_{p_2,r_2})}.
    \end{align*}
    Using Lemma \ref{Lemari¨¦-Rieusset}, we have
 \begin{align*} (-\triangle)^{-1}div(\triangle\phi\nabla\phi)\in L^1_T(B^{s_2+\frac{3}{2}-2\varepsilon_1}_{q_2,r_2}).\end{align*}
 Thus, \begin{align}\label{PP}
 P\in L^\infty_T(B^{s_1+1}_{p_1,r_1}+L^\infty_L)+L^1_T(B^{s_2+\frac{3}{2}-2\varepsilon_1}_{q_2,r_2})\hookrightarrow
 L^1_T(L^\infty_L),
  \end{align}where we have used $B^{s_1+1}_{p_1,r_1}\hookrightarrow L^\infty,$ and $B^{s_2+\frac{3}{2}-2\varepsilon_1}_{q_2,r_2}\hookrightarrow L^\infty.$

 Finally, it is easy to see that $(u,n,p,P,\phi)$ satisfies the $ENPP$ system.\\

For the uniqueness part of Theorem \ref{a1}, we first prove the uniqueness for the $\widetilde{ENPP}$ system, and then show that the solution for the $ENPP$ system also solves
 the $\widetilde{ENPP}$ system.
\subsection{Uniqueness for the $\widetilde{ENPP}$ system}
Let $\delta u=u_2-u_1,~ \delta n=n_2-n_1, \delta p=p_2-p_1,~\delta \psi=\psi_2-\psi_1.$\\
Denote $$E_i(t)\triangleq
  \|u_i\|_{\widetilde{L}^\infty_t(B^{s_1}_{p_1,r_1})}+\|n_i\|_{\widetilde{L}^\infty_t(B^{s_2}_{p_2,r_2})}+
  \|p_i\|_{\widetilde{L}^\infty_t(B^{s_2}_{p_2,r_2})},~~i=1,2,$$
  $$F_0\triangleq\|u_{02}-u_{01}\|_{B^{s'_1}_{p_1,r_1}}+\|n_{02}-n_{01}\|_{B^{s_2-1}_{p_2,r_2} }+
  \|p_{02}-p_{01}\|_{B^{s_2-1}_{p_2,r_2}},$$
  $$F(t)\triangleq \|\delta u\|_{\widetilde{L}^\infty_t(B^{s'_1}_{p_1,r_1})}+\|\delta n \|_{\widetilde{L}^\infty_t(B^{s_2-1}_{p_2,r_2}) \cap \widetilde{L}^1_t(B^{s_2+1}_{p_2,r_2})}+
  \|\delta p\|_{\widetilde{L}^\infty_t(B^{s_2-1}_{p_2,r_2}) \cap \widetilde{L}^1_t(B^{s_2+1}_{p_2,r_2})},$$
  with $s'$ denoted as in (\ref{s'}).\\
  Uniqueness for the $\widetilde{ENPP}$ system is a straightforward corollary of the following lemma.
  \begin{lemm}\label{mm}
    Let $s_1,p_1,r_1,s_2,p_2,r_2$ be as in the statement of Theorem \ref{a1}. Suppose
that we are given two solutions of the $\widetilde{ENPP}$ system
$$(u_i,n_i,p_i)\in X(T),~~i=1,2,$$
 with initial date $(u_{0i},n_{0i},p_{0i})\in B^{s_1}_{p_1,r_1}\times \big(B^{s_2}_{p_2,r_2}\big)^2.$ We then have for a.e. $t\in[0,T],$
\begin{align}\label{wending}
    F(t)\leq C(1+T)F_0exp\Big((1+T)^2C_{E_1(T)+E_2(T)}t\Big),
\end{align}	
where $C_{E_1(T)+E_2(T)}$ is a constant depending on $E_1(T)+E_2(T).$  \end{lemm}
\noindent{Proof.}
It is obvious that $\delta u$ solves
\begin{align}
     (\delta u)_t+u_2\cdot \nabla (\delta u)=-\delta u\nabla u_1-\Pi(\delta u,u_1+u_2)+\mathcal{P}\big((n_2-p_2)\delta \psi\big)
   +\mathcal{P}\big((\delta n-\delta p)\psi_1\big).
   \end{align}
According to Lemma \ref{ts}, the following inequality holds true:
\begin{align*}
  \|\delta u\|_{\widetilde{L}^\infty_t(B^{s'_1}_{p_1,r_1})}\lesssim & \|u_{02}-u_{01}\|_{B^{s'_1}_{p,r}}+C\int_0^t \|u_2\|_{B^{s_1}_{p_1,r_1}}\|\delta u\|_{\widetilde{L}_{t'}^{\infty}(B^{s'}_{p_1,r_1})}dt'+
  \|\delta u\nabla u_1\|_{\widetilde{L}^1_t(B^{s'_1}_{p_1,r_1})}\\\nonumber&+
  \|\Pi(\delta u,u_1+u_2)\|_{\widetilde{L}^1_t(B^{s'_1}_{p_1,r_1})}
    +\|\mathcal{P}\big((n_2-p_2)\delta \psi\big)\|_{\widetilde{L}^1_t(B^{s'_1}_{p_1,r_1})}
    \|\mathcal{P}\big((\delta n-\delta p)\psi_1\big)\|_{\widetilde{L}^1_t(B^{s'_1}_{p_1,r_1})}.
\end{align*}
Similar to (\ref{i}), we have,
   \begin{align*}
     \|\delta u\nabla u_1\|_{\widetilde{L}^1_t(B^{s'_1}_{p_1,r_1})}+\|\Pi(\delta u,u_1+u_2)\|_{\widetilde{L}^1_t(B^{s'_1}_{p_1,r_1})}
     \lesssim &\|\delta u\nabla u_1\|_{L^1_t(B^{s'_1}_{p_1,r_1})}+\|\Pi(\delta u,u_1+u_2)\|_{L^1_t(B^{s'_1}_{p_1,r_1})}
     \\\lesssim &\int_0^t\|\delta u\|_{B^{s'_1}_{p_1,r_1}}
     (\|u_1\|_{B^s_{p_1,r_1}}
     +\|u_2\|_{B^s_{p_1,r_1}})dt'.
   \end{align*}	
   By a similar argument as in the proof of (\ref{ii})-(\ref{iii}) , we get
   \begin{align*}
    \|\mathcal{P}\big((n_2-p_2)\delta \psi\big)\|_{\widetilde{L}^1_t(B^{s'_1}_{p_1,r_1})}\lesssim
    &\|n_2-p_2\|_{\widetilde{L}^{4}_t(B^{s_2+\frac{1}{2}}_{p_2,r_2})}
    \|\delta n-\delta p\|_{\widetilde{L}^{\frac{4}{3}}_t(B^{s_2-\frac{1}{2}}_{p_2,r_2})}\\
    \lesssim
    &\|n_2-p_2\|_{\widetilde{L}^{4}_t(B^{s_2+\frac{1}{2}}_{p_2,r_2})}
    \|\delta n-\delta p\|_{\widetilde{L}^{1}_t(B^{s_2-1}_{p_2,r_2})}^{\frac{1}{2}}\|\delta n-\delta p\|_{\widetilde{L}^{2}_t(B^{s_2}_{p_2,r_2})}^{\frac{1}{2}}\\ \lesssim
    &\delta_1^{-1}\|n_2-p_2\|_{\widetilde{L}^{4}_t(B^{s_2+\frac{1}{2}}_{p_2,r_2})}^{2}
    \|\delta n-\delta p\|_{\widetilde{L}^{1}_t(B^{s_2-1}_{p_2,r_2})}+\delta_1\|\delta n-\delta p\|_{\widetilde{L}^{2}_t(B^{s_2}_{p_2,r_2})}\\ \lesssim &\delta_1^{-1}\|n_2-p_2\|_{\widetilde{L}^{4}_T(B^{s_2+\frac{1}{2}}_{p_2,r_2})}^{2}
    \|\delta n-\delta p\|_{L^{1}_t(B^{s_2-1}_{p_2,r_2})}+\delta_1\|\delta n-\delta p\|_{\widetilde{L}^{2}_t(B^{s_2}_{p_2,r_2})},\\
     \|\mathcal{P}\big((\delta n-\delta p)\psi_1\big)\|_{\widetilde{L}^1_t(B^{s'_1}_{p_1,r_1})}\lesssim
    &\|\delta n-\delta p\|_{\widetilde{L}^{1}_t(B^{s_2+\frac{1}{2}}_{p_2,r_2})}
    \|n_1-p_1\|_{\widetilde{L}^{\infty}_t(B^{s_2-\frac{1}{2}}_{p_2,r_2})}\\
    \lesssim&
    \|\delta n-\delta p\|_{\widetilde{L}^{1}_t(B^{s_2-1}_{p_2,r_2})}^{\frac{1}{4}}\|\delta n-\delta p\|_{\widetilde{L}^{1}_t(B^{s_2+1}_{p_2,r_2})}^{\frac{3}{4}}
    \|n_1-p_1\|_{\widetilde{L}^{\infty}_t(B^{s_2}_{p_2,r_2})}\\
     \lesssim&
   \delta_1^{-3}\|\delta n-\delta p\|_{\widetilde{L}^{1}_t(B^{s_2-1}_{p_2,r_2})} \|n_1-p_1\|_{\widetilde{L}^{\infty}_t(B^{s_2}_{p_2,r_2})}^{4}+\delta_1
    \|\delta n-\delta p\|_{\widetilde{L}^{1}_t(B^{s_2+1}_{p_2,r_2})}\\\leq&
    \delta_1^{-3}\|n_1-p_1\|_{\widetilde{L}^{\infty}_T(B^{s_2}_{p_2,r_2})}^{4}\|\delta n-\delta p\|_{L^{1}_t(B^{s_2-1}_{p_2,r_2})} +\delta_1
    \|\delta n-\delta p\|_{\widetilde{L}^{1}_t(B^{s_2+1}_{p_2,r_2})}.
   \end{align*}
   Hence, we obtain
   \begin{align}\label{444}
  \|\delta u\|_{\widetilde{L}^\infty_t(B^{s'_1}_{p_1,r_1})}\lesssim ~& \|u_{02}-u_{01}\|_{B^{s'_1}_{p_1,r_1}}+
  \int_0^t[\big(\|u_1\|_{B^s_{p_1,r_1}}+\|u_2\|_{B^{s_1}_{p_1,r_1}}\big)\|\delta u\|_{\widetilde{L}^\infty_t(B^{s'_1}_{p_1,r_1})}
    \\\nonumber+&\big(\delta_1^{-1}\|n_2-p_2\|_{\widetilde{L}^{4}_T(B^{s_1+\frac{1}{2}}_{p_2,r_2})}^{2}
    +\delta_1^{-3}\|n_1-p_1\|_{\widetilde{L}^{\infty}_T(B^{s_2}_{p_2,r_2})}^{4}\big)\|\delta n-\delta p\|_{B^{s_1-1}_{p_2,r_2}}]dt'\\\nonumber&+\delta_1
    \big(\|\delta n-\delta p\|_{\widetilde{L}^{2}_t(B^{s_2}_{p_2,r_2})}+\|\delta n-\delta p\|_{\widetilde{L}^{1}_t(B^{s_2+1}_{p_2,r_2})}\big).
\end{align}

As regards $\delta n,$ note that
 \begin{align*}
     (\delta n)_t-\triangle (\delta n)&+T_{u_2}\nabla \delta n+T_{\nabla \delta n}u_2 +\nabla R(u_2,\delta n)\\&+T_{\delta u} \nabla n_1+T_{\nabla n_1}\delta u +\nabla R(\delta u, n_1)=-\nabla\cdot\big((\delta n)\psi_1+n_2 (\delta \psi)\big).
   \end{align*}
 Applying Lemma \ref{dts} yields
 \begin{align}
  &\|\delta n\|_{\widetilde{L}^\infty_t(B^{s_2-1}_{p_2,r_2})}
  +\|\delta n\|_{\widetilde{L}^1_t(B^{s_2+1}_{p_2,r_2})}\\\nonumber\lesssim&
  (1+t)\Big(\|n_{02}-n_{01}\|_{B^{s_2-1}_{p_2,r_2}}+\|\nabla\cdot\big((\delta n)\psi_1\big)\|_{\widetilde{L}^1_t(B^{s_2-1}_{p_2,r_2})}+
   \|\nabla\cdot \big(n_2 (\delta \psi)\big)\|_{\widetilde{L}^1_t(B^{s_2-1}_{p_2,r_2})}\\\nonumber&+
   \|T_{\delta u} \nabla n_1+T_{\nabla n_1}\delta u +\nabla R(\delta u, n_1)\|_{\widetilde{L}^1_t(B^{s_2-1}_{p_2,r_2})}+\|T_{u_2}\nabla \delta n+T_{\nabla \delta n}u_2 +\nabla R(u_2,\delta n)\|_{\widetilde{L}^1_t(B^{s_2-1}_{p_2,r_2})}
  .
  \end{align}
  A direct calculation similar to (\ref{44})-(\ref{77}) yields
  \begin{align}
   \|T_{u_2}\nabla \delta n+T_{\nabla \delta n}u_2 +\nabla R(u_2,\delta n)\|_{\widetilde{L}^1_t(B^{s_2-1}_{p_2,r_2})}
   \lesssim&
   \|u_2\|_{\widetilde{L}^\infty_t(B^{s'_1}_{p_1,r_1})}
   \|\delta n\|_{\widetilde{L}^1_t(B^{s_2}_{p_2,r_2})}\\\nonumber\leq&
  \delta_2^{-1}\|u_2\|_{\widetilde{L}^\infty_t(B^{s_1}_{p_1,r_1})}^{2}
   \|\delta n\|_{L^1_t(B^{s_2-1}_{p_2,r_2})}+\delta_2\|\delta n\|_{\widetilde{L}^1_t(B^{s_2+1}_{p_2,r_2})},
   \end{align}
   \begin{align}
   \|T_{\delta u} \nabla n_1+T_{\nabla n_1}\delta u +\nabla R(\delta u, n_1)
    \|_{\widetilde{L}^1_t(B^{s_2-1}_{p_2,r_2})}
\lesssim&
   \int_0^t\|\delta u \|_{B^{s'_1}_{p_1,r_1}}
   \| n_1\|_{B^{s_2}_{p_2,r_2}}dt',\end{align}
\begin{align}
  \|\nabla\cdot\big((\delta n)\psi_1\big)\|_{\widetilde{L}^1_t(B^{s_2-1}_{p_2,r_2})}\lesssim&
  \|\delta n\|_{\widetilde{L}^1_t(B^{s_2}_{p_2,r_2})}
   \|n_1-p_1\|_{\widetilde{L}^{\infty}_t(B^{s_2-\frac{1}{2}}_{p_2,r_2})}\\\nonumber\leq&
   \delta_2^{-1}\|\delta n\|_{L^1_t(B^{s_2-1}_{p_2,r_2})}
   \|n_1-p_1\|_{\widetilde{L}^{\infty}_t(B^{s_2}_{p_2,r_2})}^2+\delta_2\|\delta n\|_{\widetilde{L}^1_t(B^{s_2+1}_{p_2,r_2})},\\
  \|\nabla\cdot \big(n_2 (\delta \phi)\big)\|_{\widetilde{L}^1_t(B^{s_2-1}_{p_2,r_2})}\lesssim&
  \|n_2 \|_{\widetilde{L}^\infty_t(B^{s_2}_{p_2,r_2})}
   \|\delta n-\delta p\|_{\widetilde{L}^{1}_t(B^{s_2-\frac{1}{2}}_{p_2,r_2})}\\\nonumber \lesssim&
   \delta_2^{-\frac{1}{3}}\|\delta n-\delta p\|_{L^{1}_t(B^{s_2-1}_{p_2,r_2})}
   \|n_2 \|_{\widetilde{L}^\infty_t(B^{s_2}_{p_2,r_2})}^{\frac{4}{3}}+\delta_2\|\delta n-\delta p\|_{\widetilde{L}^{1}_t(B^{s_2+1}_{p_2,r_2})}.
 \end{align}
 Hence, we have
 \begin{align}\label{555}
  &\|\delta n\|_{\widetilde{L}^\infty_t(B^{s_2-1}_{p_2,r_2})}
  +\|\delta n\|_{\widetilde{L}^1_t(B^{s_2+1}_{p_2,r_2})}\\\nonumber\lesssim
  &(1+t)\Big(\|n_{02}-n_{01}\|_{B^{s_2-1}_{p_2,r_2}}+\int_0^t[\big(
 \delta_2^{-1}\|u_2\|_{\widetilde{L}^\infty_t(B^{s_1}_{p_1,r_1})}^{2}+\delta_2^{-1}
   \|n_1-p_1\|_{\widetilde{L}^{\infty}_t(B^{s_2}_{p_2,r_2})}^2\\\nonumber&~~+\delta_2^{-\frac{1}{3}}
   \|n_2 \|_{\widetilde{L}^\infty_t(B^{s_2}_{p_2,r_2})}^{\frac{4}{3}}\big)\big(\|\delta n\|_{B^{s_2-1}_{p_2,r_2}}+\|\delta p\|_{B^{s_2-1}_{p_2,r_2}}\big)+\| n_1\|_{B^{s_2}_{p_2,r_2}}\|\delta u \|_{B^{s'_1}_{p_1,r_1}}]
   dt'\Big)\\\nonumber&
   +(1+t)\delta_2\big(\|\delta n\|_{\widetilde{L}^{1}_t(B^{s_2+1}_{p_2,r_2})}+\|\delta p\|_{\widetilde{L}^{1}_t(B^{s_2+1}_{p_2,r_2})}\big).
  \end{align}
  Similarly,
  \begin{align}\label{666}
  &\|\delta p\|_{\widetilde{L}^\infty_t(B^{s_2-1}_{p_2,r_2})}
  +\|\delta p\|_{\widetilde{L}^1_t(B^{s_2+1}_{p_2,r_2})}\\\nonumber\lesssim
  &(1+t)\Big(\|p_{02}-p_{01}\|_{B^{s_2-1}_{p_2,r_2}}+\int_0^t[\big(
  \delta_2^{-1}\|u_2\|_{\widetilde{L}^\infty_t(B^{s_1}_{p_1,r_1})}^{2}+\delta_2^{-1}
   \|n_1-p_1\|_{\widetilde{L}^{\infty}_t(B^{s_2}_{p_2,r_2})}^2\\\nonumber&~~+\delta_2^{-\frac{1}{3}}
   \|p_2 \|_{\widetilde{L}^\infty_t(B^{s_2}_{p_2,r_2})}^{\frac{4}{3}}\big)\big(\|\delta n\|_{B^{s_2-1}_{p_2,r_2}}+\|\delta p\|_{B^{s_2-1}_{p_2,r_2}}\big)+\| p_1\|_{B^{s_2}_{p_2,r_2}}\|\delta u \|_{B^{s'_1}_{p_1,r_1}}]
   dt'\Big)\\\nonumber&
   +(1+t)\delta_2\big(\|\delta n\|_{\widetilde{L}^{1}_t(B^{s_2+1}_{p_2,r_2})}+\|\delta p\|_{\widetilde{L}^{1}_t(B^{s_2+1}_{p_2,r_2})}\big).
  \end{align}
  Choosing $\delta_1=c,~\delta_2=c(1+T)^{-1},$ and
  plugging the inequalities (\ref{555}) and (\ref{666}) into (\ref{444}), we eventually get
  \begin{align*}
    F(t)\leq C (1+T)F_0+(1+T)^2\int_0^tC_{E_1(T)+E_2(T)}F(t')dt'.
  \end{align*}
  Gronwall's lemma implies the desired result (\ref{wending}).
  \subsection{Uniqueness for the $ENPP$ system}
Subsequent the above uniqueness of the $\widetilde{ENPP}$ system, to complete the uniqueness part of Theorem \ref{a1}, it is sufficient to show that solutions to the $ENPP$ system also solves the $\widetilde{ENPP}$ system.
 \begin{lemm}
   Let $(u,n,p,P,\phi)$ satisfy the $ENPP$ system on $[0,T]\times \mathbb{R}^d.$ Assume that
\begin{align}\label{XX}
  (u,n,p)\in X(T),~and ~(P,\phi)\in Y_\alpha(T)~
  \textit{for some}~\alpha\in(0,1).
\end{align}
   Then $(u,n,p)$ satisfies the $\widetilde{ENPP}$ system, and moreover,
   $$\nabla \phi=\nabla(-\triangle)^{-1}(p-n),$$
   $$\nabla P=\Pi(u,u)+(Id-\mathcal{P})\big((n-p)\nabla(-\triangle)^{-1}(p-n)\big).$$
 \end{lemm}
 \noindent{Proof.} According to Lemma \ref{possion}, there exists a function $\phi_0\in L^\infty_T( B^{s_2+2}_{p_2,r_2}+L^\infty_{p^{\alpha_1}})$ for some $\alpha_1\in(0,1)$ satisfing
  \begin{align}
  \left\{
    \begin{array}{l}
    \triangle \phi_0=n-p, \\[1ex]
    \nabla \phi_0=\nabla(-\triangle)^{-1}(p-n).
    \end{array}
    \right.
  \end{align}
 Let $$P_{0}=P_{\pi}(u,u)-(-\triangle)^{-1}div(\triangle\phi_0\nabla\phi_0).$$
 Thanks to (\ref{PP}), $P_0\in
 L^1_T(L^\infty_L).$\\
 Note that
  \begin{align*}
    \triangle \phi=\triangle \phi_{0}=n-p.
  \end{align*}
  Hence $\phi-\phi_{0}$ is a harmonic polynomial. (\ref{XX}) and Lemma \ref{possion} guarantee that 
  $$\phi-\phi_0\in L^\infty_T(B^{s_2+2}_{p_2,r_2}+L^L_{p^{\beta}})\hookrightarrow L^\infty_T(L^L_{p^{\beta}}),$$ with $\beta=max\{\alpha,\alpha_1\}.$ \\
  This entails that $\phi-\phi_0$ depends only on t, and thus
  \begin{align}\label{pphi}
    \nabla \phi=\nabla\phi_{n-p}=\nabla(-\triangle)^{-1}(p-n).
  \end{align}
  Next applying the operator $div$ to the first equation of the $ENPP$ system, we get
  \begin{align*}
    -\triangle P=div(u\cdot\nabla u)-div~(\triangle\phi\nabla\phi)
    =-\triangle P_0.
  \end{align*}
  Note that 
  $P-P_0$ is in
  $L^1_T(L^\infty_L).$ Similar arguments as that for $\phi-\phi_0$ yield that
  \begin{align}\label{pp}
    \nabla P=\nabla P_0
    =\Pi(u,u)+(I-\mathcal{P})\big((n-p)\nabla(-\triangle)^{-1}(p-n)\big).
  \end{align}
  Thus, we conclude that  $(u,n,p)$ satisfies the $\widetilde{ENPP}$ system.
  \subsection{Properties of $n$ and $p$}
  First, note that  $$\partial_tu+u\cdot\nabla u\in \widetilde{L}^1_T(B^{s_1}_{p_1,r_1}),~\partial_tn+u\cdot\nabla n\in \widetilde{L}^1_T(B^{s_2}_{p_2,r_2}),~\partial_tp+u\cdot\nabla p\in \widetilde{L}^1_T(B^{s_2}_{p_2,r_2}).$$
  Following along the arguments in Theorem $3.19$ of \cite{keben}, we can show that \begin{align}
  u\in C([0,T];B^{s_1}_{p_1,r_1}),~ if ~r_1<\infty,~ or ~u\in C([0,T];B^{\tilde{s}_1}_{p_1,r_1}),~ if ~r_1=\infty,~\tilde{s}_1<s_1,\label{au}\\
   (n,p)\in \big(C([0,T];B^{s_2}_{p_2,r_2})\big)^2,~ if ~r_2<\infty,~or ~(n,p)\in  \big(C([0,T];B^{\tilde{s}_2}_{p_2,r_2})\big)^2,~ if ~r_2=\infty,~\tilde{s}_2<s_2.\label{auv}\end{align}
We then give a proof modeled after that of Lemma 1 in \cite{Schmuck} to show that $n,p\geq 0.$\\
  Denote $$x^+=\max\{x,0\}, ~~and~~x^-=\max\{-x,0\}.$$
 Suppose $(u,n,p,P,\phi)$ satisies the $ENPP$ system on $[0,T]\times \mathbb{R}^d.$ We introduce the following auxiliary problem
   \begin{align}\label{v}
  \left\{
    \begin{array}{l}
    \partial_tv+u\cdot\nabla v-\triangle v=-\nabla\cdot(v^+\nabla \phi), \\[1ex]
    v|_{t=0}=n_0.
    \end{array}
    \right.
  \end{align}
  We test (\ref{v}) with $(v^-)^{p_2-1}.$ After integrating by parts, we obtain
  \begin{align*}
    \frac{1}{p_2}\frac{d}{dt}\|v^-\|_{L^{p_2}}^{p_2}+(p_2-1)\int_{\mathbb{R}^d} |\nabla v^-|^2|v^-|^{p_2-2}dx=0,
  \end{align*}
  where we have used the fact
  \begin{align*}
    \big(u\cdot\nabla v,(v^-)^{p_2-1}\big)=\big(u\cdot\nabla v^+,(v^-)^{p_2-1}\big)+\big(u\cdot\nabla v^-,(v^-)^{p_2-1}\big)=0.
  \end{align*}Note that $n_0\geq 0.$ For $t\in[0,T],$ a time integration yields
  \begin{align*}
    \|v^-(t)\|_{L^{p_2}}\leq \|n_0^-\|_{L^{p_2}}=0.
  \end{align*}
  Thus, $v=v^+$ satisfies
   \begin{align}\label{vv}
  \left\{
    \begin{array}{l}
    \partial_tv+u\cdot\nabla v-\triangle v=-\nabla\cdot(v\nabla \phi), \\[1ex]
    v|_{t=0}=n_0.
    \end{array}
    \right.
  \end{align}
  Hence, $n=v\geq 0,~a.e.~on ~[0,T]\times \mathbb{R}^d.$\\
  Repeating the same steps for $p$ implies $p\geq 0,~a.e.~on ~[0,T]\times \mathbb{R}^d,$
  and this completes the proof of Theorem \ref{a1}.
  \begin{rema}
  We point out that under the conditions $(u,\nabla \phi, \triangle \phi)\in \big(L^\infty_T(L^\infty)\big)^3$ and $n_0\in L^{p_2},$ the auxiliary system (\ref{v}) (or (\ref{vv})) has a unique solution $v\in L^\infty_T(L^{p_2}),~\nabla v\in L^r_T(L^{p_2})~with~1\leq r<2.$ The proof is classical and is thus omitted.
  \end{rema}
 \section{Proof of Theorem \ref{a2}}
  This section is devoted to the proof of the continuation criterion claimed in Theorem \ref{a2}.
  Suppose that $s_1,p_1,r_1,s_2,p_2,r_2$ defined as in Theorem \ref{a1}, $s_2>\frac{d}{p_2},$ or $s_2=\frac{d}{p_2},r_2=1,$ and $q_2$ satisfies $\frac{1}{q_2}=\frac{1}{p_2}-\frac{1}{d}.$ Since $s_1> 1+\frac{d}{p_1},~s_2>\frac{d}{p_2},$ or $s_2=\frac{d}{p_2},~r_2=1,$ 
(\ref{au}) and (\ref{auv}) yield
  $$(u,n,p)\in C([0,T];L^{b_1})\times\big(C([0,T];L^a)\big)^2,$$
  with $p_1\leq b_1\leq\infty,$ and $p_2\leq a\leq\infty.$
We then have the following lemmas:
  \begin{lemm}\label{baopo}
  Assume that the ENPP system has a solution $(u,n,p,P,\phi)\in X(T)\times Y_\alpha(t),$ for some $\alpha\in(0,1).$
  Then we have for any $0\leq t< T,$
\begin{align}
\|n(t)\|_{L^a}+\|p(t)\|_{L^a}\leq 2(\|n_0\|_{L^a}+\|p_0\|_{L^a}),\label{4.1}\\
\|n(t)\|_{L^\infty}+\|p(t)\|_{L^\infty}\leq 2(\|n_0\|_{L^\infty}+\|p_0\|_{L^\infty}),\label{4.2}\\
\|\nabla \phi(t)\|_{L^b}\lesssim \|n_0\|_{L^{\frac{bd}{b+d}}}+\|p_0\|_{L^{\frac{bd}{b+d}}},\label{4.3}\\
\|\nabla \phi(t)\|_{L^\infty}\lesssim \|n_0\|_{L^{p_2}}+\|p_0\|_{L^{p_2}}+ \|n_0\|_{L^{\infty}}+\|p_0\|_{L^{\infty}},\label{4.4}
\end{align}
 where $p_2\leq a<\infty,$ and $q_2\leq b<\infty.$
\end{lemm}

 \noindent{Proof.} By multiplying both sides of the second equation of the ENPP system by $| n|^{a-2} n$ and integrating over $[0,t]\times\mathbb{R}^d,$ we get
 \begin{align}\label{n1}
   \frac{1}{a}\|n(t)\|_{L^a}^a\leq\frac{1}{a}\|n_0\|_{L^a}^a-\frac{a-1}{a}
   \int_0^t\int_{\mathbb{R}^d}\triangle \phi|n|^adxdt',
 \end{align}
 where we have used the estimates
 \begin{align*}
-\int_{\mathbb{R}^d}| n|^{a-2} n\triangle ndx=&
(a-1)\int_{\mathbb{R}^d}|n|^{a-2}|\nabla n|^2dx\geq 0,\\
\int_{\mathbb{R}^d}| n|^{a-2} n(u\cdot \nabla n)dx=&
-\frac{1}{a}\int_{\mathbb{R}^d}(div~u)|n|^adx=0,\\
-\int_{\mathbb{R}^d}| n|^{a-2} n\nabla\cdot(n\nabla\phi)dx=&
-\int_{\mathbb{R}^d}\nabla\phi\frac{1}{a}\nabla|n|^adx-\int_{\mathbb{R}^d}\triangle\phi|n|^adx\\=&
-\int_{\mathbb{R}^d}\frac{a-1}{a}\triangle \phi|n|^adx.
 \end{align*}
 Repeating the same steps for $p$ yields
\begin{align}\label{p1}
   \frac{1}{a}\|p(t)\|_{L^a}^a\leq\frac{1}{a}\|p_0\|_{L^a}^a+\frac{a-1}{a}
   \int_0^t\int_{\mathbb{R}^d}\triangle \phi|p|^adxdt'.
 \end{align}
Adding up $(\ref{n1})$ and (\ref{p1}), we get
\begin{align*}
 \frac{1}{a}(\|n(t)\|_{L^a}^a+\|p(t)\|_{L^a}^a)\leq &\frac{1}{a}(\|n_0\|_{L^a}^a+\|p_0\|_{L^a}^a)+\frac{a-1}{a}
   \int_0^t\int_{\mathbb{R}^d}\triangle \phi(|p|^a-|n|^a)dxdt'\\
   \leq &\frac{1}{a}(\|n_0\|_{L^a}^a+\|p_0\|_{L^a}^a)+\frac{a-1}{a}
   \int_0^t\int_{\mathbb{R}^d}(n-p)(p^a-n^a)dxdt'\\
   \leq &\frac{1}{a}(\|n_0\|_{L^a}^a+\|p_0\|_{L^a}^a),
\end{align*}
where we have used the non-negativity of $n,p.$
This thus leads to
\begin{align*}
\|n(t)\|_{L^a}+\|p(t)\|_{L^a}\leq & 2^{1-\frac{1}{a}}(\|n(t)\|_{L^a}^a+\|p(t)\|_{L^a}^a)^{\frac{1}{a}}
\\\leq & 2^{1-\frac{1}{a}}(\|n_0\|_{L^a}^a+\|p_0\|_{L^a}^a)^{\frac{1}{a}}\\\leq & 2(\|n_0\|_{L^a}+\|p_0\|_{L^a}).
\end{align*}
Hence, (\ref{4.1}) holds.
Then passing to the limit as $a$ tends to infinity gives the inequality $(\ref{4.2})$,
and the inequality $(\ref{4.3})$ is just an application of Lemma \ref{Lemari¨¦-Rieusset} and the inequality $(\ref{4.1})$.
In order to prove $(\ref{4.4}),$ we split $\triangle \phi$ into low and high frequencies \begin{align*}
\nabla \phi=&\triangle_{-1}\nabla \phi+(Id-\triangle_{-1})\nabla \phi
.\end{align*}
By virtue of Lemmas \ref{Bi} and \ref{Lemari¨¦-Rieusset}, we deduce that
\begin{align*}
\|(\triangle_{-1}\nabla \phi)(t)\|_{L^\infty}\lesssim\|(\triangle_{-1}\nabla \phi)(t)\|_{L^{q_2}}\lesssim\|\big(\triangle_{-1}(p-n)\big)(t)\|_{L^{p_2}}\lesssim\|(p-n)(t)\|_{L^{p_2}}.
\end{align*}
As $\nabla\nabla(-\triangle)^{-1}$ is a Calderon-Zygmund operator and $p_2<d,$ we have
\begin{align*}
\|\big(Id-\triangle_{-1})\nabla \phi\big)(t)\|_{L^\infty}\lesssim &
\|\mathcal{F}^{-1}\big(\chi(0)-\chi(\xi)\big)\ast\nabla \phi(t)\|_{L^\infty}
\\\lesssim & \|\int_0^1\frac{1}{\tau^d}h(\frac{\cdot}{\tau})d\tau\ast\nabla(\nabla \phi)(t)\|_{L^\infty}\\
\lesssim &\|\int_0^1\frac{1}{\tau^d}h(\frac{\cdot}{\tau})d\tau\|_{L^{l'}}\|\nabla\big(\nabla \triangle^{-1}(p-n)\big)(t)\|_{L^{l}}\\
\lesssim &\int_0^1\frac{1}{\tau^{\frac{d}{l}}}d\tau\|(n-p)(t)\|_{L^{l}}\\
\lesssim &\|(n-p)(t)\|_{L^{p_2}}+\|(n-p)(t)\|_{L^\infty}\\\lesssim &\|n_0\|_{L^{p_2}}+\|p_0\|_{L^{p_2}}+\|n_0\|_{L^\infty}+\|p_0\|_{L^\infty},
\end{align*}
 where $\mathcal{F}(h)=\chi'(\xi),$ $d< l<\infty,$ and $l'$ is the conjugate of $l.$
Thus the inequality $(\ref{4.4})$ is proved. \qed
 \begin{lemm}\label{baopo2}
Under the assumption of Lemma \ref{baopo}, for any $b_1>d,~b_1\geq p_1,$ if
 \begin{align}
   \|u\|_{L^\infty_T(L^{b_1})}<\infty,
 \end{align}
  then we have 
 \begin{align}\label{qq}
   \|\nabla n\|_{L^1_T(L^{a_1})}+
   \|\nabla p\|_{L^1_T(L^{a_1})}\leq C_0(T)<\infty,
 \end{align}
 where $p_2\leq a_1\leq\infty,$ and $C_0(T)$ depends on $\|u_0\|_{B^{s_1}_{p_1,r_1}},~\|n_0\|_{B^{s_2}_{p_2,r_2}},~
  \|p_0\|_{B^{s_2}_{p_2,r_2}}$, $T$ and $\|u\|_{L^\infty_T(L^{b_1})}.$~~~
 \end {lemm}
\noindent{Proof.} Note that
\begin{align*}
  \nabla n =\nabla e^{t\triangle}n_0-\int_0^te^{(t-t')\triangle}\nabla (u\cdot \nabla n+\nabla n\nabla \phi+n\triangle\phi)dt'.
\end{align*}
It is easy to obtain
\begin{align*}
  &\|(\nabla n)(\tau)\|_{L^{a_1}}\\\lesssim&\| \mathcal{F}^{-1}(e^{-\tau|\xi|^2}\xi)\|_{L^1}\|n_0\|_{L^{a_1}}+ \int_0^\tau\|\mathcal{F}^{-1}(e^{(t'-\tau)|\xi|^2}\xi)\big\|_{L^{{b_1}'}}\|u(t')\|_{L^{b_1}}\| \nabla n(t')\|_{L^{a_1}}dt'\\
  &+\int_0^\tau\big\|\mathcal{F}^{-1}(e^{(t'-\tau)|\xi|^2}\xi)\big\|_{L^1}\big(\|\nabla \phi\|_{L^\infty}\|\nabla n\|_{L^{a_1}}+\|\triangle \phi\|_{L^\infty}\|n\|_{L^{a_1}}\big)(t')dt'\\\lesssim& \tau^{-\frac{1}{2}}\|n_0\|_{L^{a_1}}+ \int_0^\tau\frac{1}{{(\tau-t')}^{\frac{1}{2}+\frac{d}{2b_1}}}\|u(t')\|_{L^{b_1}}\| \nabla n(t')\|_{L^{a_1}}dt'\\
  &+\int_0^\tau\frac{1}{{(\tau-t')}^{\frac{1}{2}}}\big(\|\nabla \phi\|_{L^\infty}\|\nabla n\|_{L^{a_1}}+\|\triangle \phi\|_{L^\infty}\|n\|_{L^{a_1}}\big)(t')dt'\\
  \lesssim&\tau^{-\frac{1}{2}}\| n_0\|_{L^{a_1}}+\tau^{\frac{1}{2}}\|\triangle \phi\|_{L^\infty_T(L^\infty)}\|n\|_{L^\infty_T(L^{a_1})}\\&+\int_0^\tau\big(\frac{1}{{(\tau-t')}^{\frac{1}{2}+\frac{d}{2b_1}}}
  \|u(t')\|_{L^{b_1}}+\frac{1}{{(\tau-t')}^{\frac{1}{2}}}\|\nabla \phi(t')\|_{L^\infty}\big)\|\nabla n(t')\|_{L^{a_1}}dt'\\
  \lesssim&\tau^{-\frac{1}{2}}\| n_0\|_{L^{a_1}}+\tau^{\frac{1}{2}}\|\triangle \phi\|_{L^\infty_T(L^\infty)}\|n\|_{L^\infty_T(L^{a_1})}\\&
  +\int_0^\tau\big(\delta_1\frac{1}{{(\tau-t')}^{(\frac{1}{2}+\frac{d}{2b_1})\gamma_1}}
  +\delta_2\frac{1}{{(\tau-t')}^{\frac{1}{2}\gamma_2}}\big)\|\nabla n(t')\|_{L^{a_1}}dt'\\&+\int_0^\tau\big(C_{\delta_1}
  \|u(t')\|^{\gamma'_1}_{L^{b_1}}+C_{\delta_2}\|\nabla \phi(t')\|^{\gamma'_2}_{L^\infty}\big)\|\nabla n(t')\|_{L^{a_1}}dt',
\end{align*}
 where $\gamma_1>1,~\gamma_2>1$ are chosen to satisfy
 \begin{align*}
   (\frac{1}{2}+\frac{d}{2b_1})\gamma_1<1,~~\frac{1}{2}\gamma_2<1,
 \end{align*}
  $\delta_1,~\delta_2$ sufficiently small will be defined later, and $x'$ denotes the conjugate of $x.$  
By means of the Young inequality for the time integral, we obtain,
\begin{align*}
  \|\nabla n\|_{L^1_t(L^{a_1})}\lesssim &t^{\frac{1}{2}}\| n_0\|_{L^{a_1}}+t^{\frac{3}{2}}\|\triangle \phi\|_{L^\infty_T(L^\infty)}\|n\|_{L^\infty_T(L^{a_1})}\\&
  +\big(\delta_1\frac{1}{1-(\frac{1}{2}+\frac{d}{2b_1})\gamma_1}t^{1-(\frac{1}{2}+\frac{d}{2b_1})\gamma_1}
  +\delta_2\frac{1}{1-\frac{1}{2}\gamma_2}t^{1-\frac{1}{2}\gamma_2}\big)\|\nabla n\|_{L^1_t(L^{a_1})}\\&+\int_0^t\int_0^\tau\big(C_{\delta_1}
  \|u(t')\|^{\gamma'_1}_{L^{b_1}}+C_{\delta_2}\|\nabla \phi(t')\|^{\gamma'_2}_{L^\infty}\big)\|\nabla n(t')\|_{L^{a_1}}dt'd\tau.
\end{align*}
Choosing $\delta_1\frac{1}{1-(\frac{1}{2}+\frac{d}{2b_1})\gamma_1}T^{1-(\frac{1}{2}+\frac{d}{2b_1})\gamma_1}
  +\delta_2\frac{1}{1-\frac{1}{2}\gamma_2}T^{1-\frac{1}{2}\gamma_2}=c\frac{1}{2}$ yields
\begin{align*}
   \|\nabla n\|_{L^1_t(L^{a_1})}\lesssim &T^{\frac{1}{2}}\| n_0\|_{L^{a_1}}+T^{\frac{3}{2}}\|\triangle \phi\|_{L^\infty_T(L^\infty)}\|n\|_{L^\infty_T(L^{a_1})}\\&
  +\big(C_{\delta_1}
  \|u(t')\|^{\gamma'_1}_{L^\infty_T(L^{b_1})}+C_{\delta_2}\|\nabla \phi(t')\|^{\gamma'_2}_{L^\infty_T(L^\infty)}\big)\int_0^t\int_0^\tau\|\nabla n(t')\|_{L^{a_1}}dt'd\tau.
   \end{align*}
Gronwall's lemma thus implies that
\begin{align*}
\|\nabla n\|_{L^1_t(L^{a_1})}\leq &C\Big(T^{\frac{1}{2}}\| n_0\|_{L^{a_1}}+T^{\frac{3}{2}}\|\triangle \phi\|_{L^\infty_T(L^\infty)}\|n\|_{L^\infty_T(L^{a_1})}\Big)\\ &\times exp\Big(\big(C_{\delta_1}
  \|u(t')\|^{\gamma'_1}_{L^\infty_T(L^{b_1})}+C_{\delta_2}\|\nabla \phi(t')\|^{\gamma'_2}_{L^\infty_T(L^\infty)}\big)t\Big).
  \end{align*}
 Similar arguments for $p,$ combining with Lemma \ref{baopo} yield the inequality $(\ref{qq})$ holds true. \qed

\begin{lemm}\label{baopo3}
  Under the assumption of Lemma \ref{baopo}, if $1<p_1<\infty,$ and  $\int_0^T\|\nabla u\|_{L^\infty}dt'<\infty,$
  then we have for any $0\leq t<T$ and $p_1\leq b_1<\infty,$
   \begin{align}\label{4.7}
   \|u\|_{L^\infty_t(L^{b_1})}&\leq C_1(T)
    exp\big(C\int_0^t\|\nabla u\|_{L^\infty}dt'\big)<\infty,\\
    \|u\|_{L^1_t(L^{\infty})}&\lesssim TC_1(T)
    exp\big(C\int_0^t\|\nabla u\|_{L^\infty}dt'\big)+\int_0^t\|\nabla u\|_{L^\infty}<\infty,
   \end{align}
   with $C_1(T)$ depending on $\|u_0\|_{B^{s_1}_{p_1,r_1}},~\|n_0\|_{B^{s_2}_{p_2,r_2}},~
  \|p_0\|_{B^{s_2}_{p_2,r_2}}$ and $T.$
   \end{lemm}
  \noindent{Proof.} We first deduce from the condition (\ref{1.9}) that there exist $l_1$ and $l_2,$ such that $p_2\leq l_1\leq \infty,~q_2\leq l_2\leq \infty,$ and $\frac{1}{l_1}+\frac{1}{l_2}=\frac{1}{p_1}.$\\ As $u$ satisfies the first equation of the $\widetilde{ENPP}$ system, we have
  \begin{align*}
    \|u\|_{L^{b_1}}\leq &\|u_0\|_{L^{b_1}}+\int_0^t\|\Pi(u,u)\|_{L^{b_1}}+\|\mathcal{P}(\triangle\phi\nabla\phi)\|_{L^{b_1}}dt'\\
    \nonumber \leq&\|u_0\|_{L^{b_1}}+C\big(\int_0^t\|u\|_{L^{b_1}}\|\nabla u\|_{L^\infty}+\|\triangle\phi\|_{L^{l_1}}\|\nabla\phi\|_{L^{l_2}}dt'\big),
  \end{align*}
  where we have used Lemma \ref{bao}.
 Lemma \ref{baopo} and Gronwall's lemma thus imply that
  \begin{align*}
    \|u\|_{L^\infty_t(L^{b_1})}\leq &\big(\|u_0\|_{L^{b_1}}+C\int_0^{t}\|\triangle\phi\|_{L^{l_1}}\|\nabla\phi\|_{L^{l_2}}dt'\big)
    exp\big(C\int_0^t\|\nabla u\|_{L^\infty}dt'\big)<\infty\\
    \leq&\big(\|u_0\|_{B^{s_1}_{p_1,r_1}}+C_{(\|n_0\|_{B^{s_2}_{p_2,r_2}}+\|p_0\|_{B^{s_2}_{p_2,r_2}})}t\big)
    exp\big(C\int_0^t\|\nabla u\|_{L^\infty}dt'\big),
  \end{align*}
  where we have used the facts $B^{s_1}_{p_1,r_1}\hookrightarrow L^{b_1},$ $B^{s_2}_{p_2,r_2}\hookrightarrow L^{l_1},$ and  $B^{s_2}_{p_2,r_2}\hookrightarrow L^{l_2}.$
 Then, the Gronwall Lemma gives the inequality (\ref{4.7}). \\
  Next, by splitting $u$ into low and high frequencies, we see that
  \begin{align*}
    \|u\|_{L^\infty}\lesssim\|\triangle_{-1}u\|_{L^\infty}+\|(I-\triangle_{-1})u\|_{L^\infty}
    \lesssim\|u\|_{L^{p_1}}+\|\nabla u\|_{L^\infty}.
  \end{align*}
  Applying the inequality (\ref{4.7}), we complete the proof of the lemma. \qed

We now turn to the proof of Theorem \ref{a2}. Note that $(u,n,p)$ also satisfies the $\widetilde{ENPP}$ system. Applying $\triangle_j$ to the first equation of the $\widetilde{ENPP}$ system yields that
\begin{align*}
(\partial_t+u\cdot \nabla )\triangle_ju+\triangle_j\Pi(u,u)=\triangle_j\mathcal{P}(\triangle\phi\nabla\phi)+R_{j1},
\end{align*}
with $R_{j1}=u\cdot \nabla\triangle_ju-\triangle_j(u\cdot \nabla)u\triangleq[u\cdot\nabla,\triangle_j]u.$\\
Using the fact that $div~u=0,$ we readily obtain
\begin{align*}
  \|\triangle_ju(t)\|_{L^{p_1}}\leq\|\triangle_ju_0\|_{L^{p_1}}+\int_0^t\|\triangle_j\Pi(u,u)\|_{L^{p_1}}
  +\|\triangle_j\mathcal{P}(\triangle\phi\nabla\phi)\|_{L^{p_1}}+\|R_{j1}\|_{L^{p_1}}dt',~j\geq-1.
\end{align*}
Multiplying both sides of the above inequality by $2^{js_1}$, taking the $l^{r_1}$ norm and using the Minkowski inequality, we obtain
\begin{align}\label{uu'}
  \|u\|_{\widetilde{L}^\infty_t(B^{s_1}_{p_1,r_1})}\lesssim
  \|u_0\|_{B^{s_1}_{p_1,r_1}}+\|\Pi(u,u)\|_{L^1_t(B^{s_1}_{p_1,r_1})}+
  \|\mathcal{P}(\triangle\phi\nabla\phi)\|_{\widetilde{L}^1_t(B^{s_1}_{p_1,r_1})}
  +\Big\|2^{js_1}\|R_{j1}\|_{L^1_t(L^{p_1})}\Big\|_{l^{r_1}}.
\end{align}

Due to Lemma \ref{jiaohuan}, we get
\begin{align}\label{RR}
  \Big\|2^{js_1}\|R_{j1}\|_{L^{p_1}}\Big\|_{l^{r_1}}\lesssim &\|\nabla u\|_{L^\infty}\|u\|_{B^{s_1}_{p_1,r_1}}.
\end{align}
By virtue of Lemma \ref{pi}, we have
\begin{align}\label{pipi}
  \|\Pi(u,u)\|_{B^{s_1}_{p_1,r_1}}\lesssim \|v\|_{C^{0,1}}\|u\|_{B^{s_1}_{p_1,r_1}}.
\end{align}
We now focus on the term $\mathcal{P}(\triangle\phi\nabla\phi).$\\
If $p_1<p_2,$ let $\frac{1}{p_3}=\frac{1}{p_1}-\frac{1}{p_2}.$
 The condition (1.9) implies that $p_2\leq q_2\leq p_3
.$ We have
\begin{align*}
  \|\mathcal{P}(\triangle\phi\nabla\phi)\|_{B^{s_1}_{p_1,r_1}}\lesssim&
  \|\triangle\phi\nabla\phi\|_{B^{s_1}_{p_1,r_1}}\\ \lesssim&
  \|T_{\triangle\phi}(Id-\triangle_{-1})\nabla\phi\|_{B^{s_1}_{p_1,r_1}}
  +\|T_{\nabla\phi}\triangle\phi+R(\triangle\phi,\nabla\phi)\|_{B^{s_1}_{p_1,r_1}}\\ \lesssim&\|\triangle\phi\|_{L^{p_3}}\|(Id-\triangle_{-1})\nabla\phi\|_{B^{s_1}_{p_2,r_1}}
  +\|\nabla\phi\|_{L^{p_3}}\|\triangle\phi\|_{B^{s_1}_{p_2,r_1}}\\ \lesssim&\|\triangle\phi\|_{L^{p_3}}\|\triangle\phi\|_{B^{s_2+\frac{1}{2}}_{p_2,r_2}}
  +\|\nabla\phi\|_{L^{p_3}}\|\triangle\phi\|_{B^{s_2+\frac{3}{2}}_{p_2,r_2}}.
  \end{align*}
  If $p_1\geq p_2,$ we have
  \begin{align*}
  \|\mathcal{P}(\triangle\phi\nabla\phi)\|_{B^{s_1}_{p_1,r_1}}\lesssim&
  \|\mathcal{P}(\triangle\phi\nabla\phi)\|_{B^{s_1-\frac{d}{p_1}+\frac{d}{p_2}}_{p_2,r_1}} \lesssim
  \|\triangle\phi\nabla\phi\|_{B^{s_2+\frac{3}{2}}_{p_2,r_2}}\\
  \lesssim&
  \|\triangle\phi\|_{L^{\infty}}\|(Id-\triangle_{-1})\nabla\phi\|_{B^{s_2+\frac{3}{2}}_{p_2,r_2}}
  +\|\nabla\phi\|_{L^{\infty}}\|\triangle\phi\|_{B^{s_2+\frac{3}{2}}_{p_2,r_2}}\\ \lesssim&\|\triangle\phi\|_{L^{\infty}}\|\triangle\phi\|_{B^{s_2+\frac{1}{2}}_{p_2,r_2}}
  +\|\nabla\phi\|_{L^{\infty}}\|\triangle\phi\|_{B^{s_2+\frac{3}{2}}_{p_2,r_2}}.
  \end{align*}
Let $a_1=\frac{p_1p_2}{p_2-p_1},
 ~if ~p_1<p_2,$ and $a_1=\infty,~if~ p_1\geq p_2.$ By using interpolation and the Young inequality, we infer that
 \begin{align}\label{phiphi}
   &\|\mathcal{P}(\triangle\phi\nabla\phi)\|_{\widetilde{L}^1_t(B^{s_1}_{p_1,r_1})}\\\nonumber \lesssim&
   \|\triangle\phi\|_{L^{\infty}_t(L^{a_1})}
   \|\triangle\phi\|_{\widetilde{L}^1_t(B^{s_2+\frac{1}{2}}_{p_2,r_2})}
  +\|\nabla\phi\|_{L^{\infty}_t(L^{a_1})}
  \|\triangle\phi\|_{\widetilde{L}^1_t(B^{s_2+\frac{3}{2}}_{p_2,r_2})}\\\nonumber  \lesssim&
  \|\triangle\phi\|_{L^{\infty}_t(L^{a_1})}
   \|\triangle\phi\|_{\widetilde{L}^1_t(B^{s_2}_{p_2,r_2})}^{\frac{3}{4}}
   \|\triangle\phi\|_{\widetilde{L}^1_t(B^{s_2+2}_{p_2,r_2})}^{\frac{1}{4}}
  +\|\nabla\phi\|_{L^{\infty}(L^{a_1})}
  \|\triangle\phi\|_{\widetilde{L}^1_t(B^{s_2}_{p_2,r_2})}^{\frac{1}{4}}
  \|\triangle\phi\|_{\widetilde{L}^1_t(B^{s_2+2}_{p_2,r_2})}^{\frac{3}{4}}\\ \nonumber
  \leq &C_{\sigma_1}\Big(\|\triangle\phi\|_{L^{\infty}_t(L^{a_1})}^{\frac{4}{3}}
   \|\triangle\phi\|_{\widetilde{L}^1(B^{s_2}_{p_2,r_2})}+
   \|\nabla\phi\|_{L^{\infty}_t(L^{a_1})}^{4}
  \|\triangle\phi\|_{\widetilde{L}^1_t(B^{s_2}_{p_2,r_2})}\Big)+\sigma_1
  \|\triangle\phi\|_{\widetilde{L}^1_t(B^{s_2+2}_{p_2,r_2})}\\ \nonumber
  \leq &C_{\sigma_1}\Big(\|\triangle\phi\|_{L^{\infty}_t(L^{a_1})}^{\frac{4}{3}}
   +
   \|\nabla\phi\|_{L^{\infty}_t(L^{a_1})}^{4}\Big)
  \|\triangle\phi\|_{L^1_t(B^{s_2}_{p_2,r_2})}+\sigma_1
  \|\triangle\phi\|_{\widetilde{L}^1_t(B^{s_2+2}_{p_2,r_2})}.
 \end{align}
 Plugging the inequalities (\ref{RR})-(\ref{phiphi}) into (\ref{uu'}), we eventually get
 \begin{align}\label{345}
 \|u\|_{\widetilde{L}^\infty_t(B^{s_1}_{p_1,r_1})}\leq &\|u_0\|_{B^{s_1}_{p_1,r_1}}+\int_0^t[C\|u\|_{C^{0,1}}
 \|u\|_{B^{s_1}_{p_1,r_1}}\\\nonumber&+C_{\sigma_1}\Big(\|\triangle\phi\|_{\widetilde{L}^{\infty}_t(L^{a_1})}^{\frac{4}{3}}+
   \|\nabla\phi\|_{\widetilde{L}^{\infty}_t(L^{a_1})}^{4}\Big)
  \|\triangle\phi\|_{B^{s_2}_{p_2,r_2}}]dt'+\sigma_1
  \|\triangle\phi\|_{\widetilde{L}^1_t(B^{s_2+2}_{p_2,r_2})}.
  \end{align}
 Similarly, applying $\triangle_j$ to the second equation of the $ENPP$ system yields that
\begin{align*}
(\partial_t+u\cdot \nabla -\triangle) \triangle_jn=-\triangle_j\nabla\cdot(n\nabla\phi)+R_{j2},
\end{align*}
with $R_{j2}=[u\cdot\nabla,\triangle_j]n.$\\
Thanks to Lemma \ref{DL}, and using the fact that $div~u=0,$ we readily obtain
\begin{align*}
  \|\triangle_jn(t)\|_{L^{p_2}}\leq\|\triangle_jn_0\|_{L^{p_2}}+\int_0^t
  \|\triangle_j\nabla\cdot(n\nabla\phi)\|_{L^{p_2}}
  +\|R_{j2}\|_{L^{p_2}}dt',~j=-1,\\
  \|\triangle_jn(t)\|_{L^{p_2}}+C\int_0^t2^{2j}\|\triangle_jn\|_{L^{p_2}}dt'\leq\|\triangle_jn_0\|_{L^{p_2}}+\int_0^t
  \|\triangle_j\nabla\cdot(n\nabla\phi)\|_{L^{p_2}}
  +\|R_{j2}\|_{L^p}dt',~j\geq0,
\end{align*}
from which it follows that
\begin{align}\label{xingxing}
  \|\triangle_jn(t)\|_{L^{p_2}}+\int_0^t2^{2j}\|\triangle_jn\|_{L^{p_2}}dt'\lesssim (1+t)\big(&\|\triangle_jn_0\|_{L^{p_2}}+\int_0^t
  \|\triangle_j\nabla\cdot(n\nabla\phi)\|_{L^{p_2}}
  \\\nonumber+&\|R_{j2}\|_{L^p}dt'\big),~j\geq-1.
\end{align}
Hence multiplying both sides of the above inequality by $2^{js_2}$ and taking the $l^{r_2}$ norm, we obtain
\begin{align*}
  &\|n\|_{\widetilde{L}^\infty_t(B^{s_2}_{p_2,r_2})}+ \|n\|_{\widetilde{L}^1_t(B^{s_2+2}_{p_2,r_2})}\\\lesssim& (1+t)\big(\|n_0\|_{B^{s_2}_{p_2,r_2}}+
  \|\nabla\cdot(n\nabla\phi)\|_{\widetilde{L}^1_t(B^{s_2}_{p_2,r_2})}
  +\Big\|2^{js_2}\|R_{j2}\|_{L^1_t(L^{p_2})}\Big\|_{l^{r_2}}\big).
\end{align*}
In view of Lemma \ref{jiaohuan}, we get
\begin{align}\label{R2}
   \Big\|2^{js_2}\|R_{j2}\|_{L^{p_2}}\Big\|_{l^{r_2}}\lesssim &\Big(\|\nabla u\|_{L^\infty}\|n\|_{B^{s_2}_{p_2,r_2}}+\|\nabla n\|_{L^a}\|\nabla u\|_{B^{s_2-1}_{b,r_2}}\Big)\\\nonumber
   \lesssim &\Big(\|\nabla u\|_{L^\infty}\|n\|_{B^{s_2}_{p_2,r_2}}+\|\nabla n\|_{L^a}\|\nabla u\|_{B^{s_2-1+\frac{d}{p_1}-\frac{d}{b}}_{p_1,r_2}}\Big)
   \\\nonumber
   \lesssim &\Big(\|\nabla u\|_{L^\infty}\|n\|_{B^{s_2}_{p_2,r_2}}+\|\nabla n\|_{L^a}\|\nabla u\|_{B^{s_1-1}_{p_1,r_1}}\Big),
\end{align}
where $a=\infty,~b=p_2,~if~p_1\leq p_2,$ and $a=\frac{p_1p_2}{p_1-p_2},~b=p_1,~if~p_1>p_2.$\\
According to Lemmas \ref{T}-\ref{R}, we have
\begin{align}\label{n2}
  \|\nabla n\nabla \phi\|_{B^{s_2}_{p_2,r_2}}\lesssim&\|\nabla n\|_{L^d}\|\nabla \phi\|_{B^{s_2}_{q_2,r_2}}+\|\nabla \phi\|_{L^\infty}\|\nabla n\|_{B^{s_2}_{p_2,r_2}}\\\nonumber
  \lesssim&\|\nabla n\|_{L^d}\|\triangle \phi\|_{B^{s_2}_{p_2,r_2}}+\|\nabla \phi\|_{L^\infty}\|\nabla n\|_{B^{s_2}_{p_2,r_2}},\\
 \label{n3}
  \|n\triangle\phi\|_{B^{s_2}_{p_2,r_2}}\lesssim&\|n\|_{L^\infty}\|\triangle\phi
  \|_{B^{s_2}_{p_2,r_2}}+\|\triangle\phi\|
  _{L^\infty}\| n\|_{B^{s_2}_{p_2,r_2}}.
\end{align}
Inserting the inequalities (\ref{R2})-(\ref{n3}) into (\ref{xingxing}) and using the Minkowski inequality, we finally get
\begin{align}\label{567}
 &\|n\|_{\widetilde{L}^\infty_t(B^{s_2}_{p_2,r_2})}+\|n\|_{\widetilde{L}^1_t(B^{s_2+2}_{p_2,r_2})}\\\nonumber
\lesssim&(1+t)\Big(\|n_0\|_{B^{s_2}_{p_2,r_2}}+\int_0^t[(\|\nabla u\|_{L^\infty}+\|\triangle\phi\|
  _{L^\infty})\|n\|_{B^{s_2}_{p_2,r_2}}+\|\nabla n\|_{L^a}\| u\|_{B^{s_1}_{p_1,r_1}} \\\nonumber &+(\|\nabla n\|_{L^d}+\|n\|_{L^\infty})\|\triangle \phi\|_{B^{s_2}_{p_2,r_2}}]dt'+
  \|\nabla \phi\|_{L^\infty_t(L^\infty)}\|\nabla n\|_{\widetilde{L}^1_t(B^{s_2}_{p_2,r_2})}\\\nonumber
\lesssim&(1+t)\Big(\|n_0\|_{B^{s_2}_{p_2,r_2}}+\int_0^t[(\|\nabla u\|_{L^\infty}+\|\triangle\phi\|
  _{L^\infty}+C_{\sigma_2}\|\nabla \phi\|_{L^\infty_t(L^\infty)}^2)\|n\|_{B^{s_2}_{p_2,r_2}} +\|\nabla n\|_{L^a}\| u\|_{B^{s_1}_{p_1,r_1}} \\\nonumber&+(\|\nabla n\|_{L^d}+\|n\|_{L^\infty})\|\triangle \phi\|_{B^{s_2}_{p_2,r_2}}]dt'+
  \sigma_2\| n\|_{\widetilde{L}^1_t(B^{s_2+2}_{p_2,r_2})}\Big),
  \end{align}
  where we have used that
  \begin{align*}
    \|\nabla \phi\|_{L^\infty_t(L^\infty)}\|\nabla n\|_{\widetilde{L}^1_t(B^{s_2}_{p_2,r_2})}\lesssim&
    \|\nabla \phi\|_{L^\infty_t(L^\infty)}\| n\|_{\widetilde{L}^1_t(B^{s_2}_{p_2,r_2})}^{\frac{1}{2}}\| n\|_{\widetilde{L}^1_t(B^{s_2+2}_{p_2,r_2})}^{\frac{1}{2}}\\
    \lesssim&C_{\sigma_2}\|\nabla \phi\|_{L^\infty_t(L^\infty)}^2\| n\|_{\widetilde{L}^1_t(B^{s_2}_{p_2,r_2})}+\sigma_2\| n\|_{\widetilde{L}^1_t(B^{s_2+2}_{p_2,r_2})}\\
    \lesssim&\int_0^tC_{\sigma_2}\|\nabla \phi\|_{L^\infty_t(L^\infty)}^2\| n\|_{B^{s_2}_{p_2,r_2}}dt'+\sigma_2\| n\|_{\widetilde{L}^1_t(B^{s_2+2}_{p_2,r_2})}.
  \end{align*}
  Similar arguments as above yield \begin{align}\label{568}
  &\|p\|_{\widetilde{L}^\infty_t(B^{s_2}_{p_2,r_2})}+\|p\|_{\widetilde{L}^1_t(B^{s_2+2}_{p_2,r_2})}\\\nonumber
\lesssim&(1+t)\Big(\|p_0\|_{B^{s_2}_{p_2,r_2}}+\int_0^t[(\|\nabla u\|_{L^\infty}+\|\triangle\phi\|
  _{L^\infty}+C_{\sigma_2}\|\nabla \phi\|_{L^\infty_t(L^\infty)}^2)\|p\|_{B^{s_2}_{p_2,r_2}} +\|\nabla p\|_{L^a}\| u\|_{B^{s_1}_{p_1,r_1}} \\\nonumber&+(\|\nabla p\|_{L^d}+\|p\|_{L^\infty})\|\triangle \phi\|_{B^{s_2}_{p_2,r_2}}]dt'+
  \sigma_2\|p\|_{\widetilde{L}^1_t(B^{s_2+2}_{p_2,r_2})}\Big).
  \end{align}
  Combining (\ref{567}), (\ref{568}) and (\ref{345}), we get
  \begin{align}
    &\|u\|_{\widetilde{L}^\infty_t(B^{s_1}_{p_1,r_1})}+\|n\|_{\widetilde{L}^\infty_t(B^{s_2}_{p_2,r_2})\cap \widetilde{L}^1_t(B^{s_2+2}_{p_2,r_2})}+\|p\|_{\widetilde{L}^\infty_t(B^{s_2}_{p_2,r_2})\cap \widetilde{L}^1_t(B^{s_2+2}_{p_2,r_2})}\\
    \nonumber\lesssim
    &(1+t)\big(\|u_0\|_{B^{s_1}_{p_1,r_1}}+\|n_0\|_{B^{s_2}_{p_2,r_2}}+\|p_0\|_{B^{s_2}_{p_2,r_2}}\big)+(1+t)\int_0^t\Big(
    \|u\|_{C^{0,1}}+\| n\|_{L^\infty}+\| p\|_{L^\infty}\\
    \nonumber&+\|\nabla n\|_{L^a}+\|\nabla p\|_{L^a}+\|\nabla n\|_{L^d}+\|\nabla p\|_{L^d}+C_{\sigma_1}\|\triangle\phi\|_{L^{\infty}_{T}(L^{a_1})}^{\frac{4}{3}}+
   C_{\sigma_1}\|\nabla\phi\|_{L^{\infty}_{T}(L^{a_2})}^{4}\\\nonumber&+\sigma_2\|\nabla \phi\|_{L^\infty_{T}(L^\infty)}^2\Big)
   \times\Big(\|u\|_{B^{s_1}_{p_1,r_1}}+\|n\|_{B^{s_2}_{p_2,r_2}}+\|p\|_{B^{s_2}_{p_2,r_2}}\Big)dt'
   \\\nonumber&+\big(\sigma_1+(1+t)\sigma_2\big)\big(\|n\|_{\widetilde{L}^1_t(B^{s_2+2}_{p_2,r_2})}+\|p\|_{\widetilde{L}^1_t(B^{s_2+2}_{p_2,r_2})}\big).
  \end{align}
  Choose $\sigma_1=c,~\sigma_2=c(1+T)^{-1}.$
  Lemmas \ref{baopo}-\ref{baopo2} imply that
  if (\ref{b1}) holds, then
  \begin{align*}
    &\|u\|_{C^{0,1}}+\| n\|_{L^\infty}+\| p\|_{L^\infty}
+\|\nabla n\|_{L^a}+\|\nabla p\|_{L^a}+\|\nabla n\|_{L^d}+\|\nabla p\|_{L^d}\\+&C_{\sigma_1}\|\triangle\phi\|_{L^{\infty}_{T}(L^{a_1})}^{\frac{4}{3}}+
   C_{\sigma_1}\|\nabla\phi\|_{L^{\infty}_{T}(L^{a_2})}^{4}+\sigma_2\|\nabla \phi\|_{L^\infty_{T}(L^\infty)}^2\Big)\in L^1_T,
  \end{align*}
Applying Gronwall's lemma thus leads to
  \begin{align*}
    \|u\|_{\widetilde{L}^\infty_T(B^{s_1}_{p_1,r_1})}+\|n\|_{\widetilde{L}^\infty_T(B^{s_2}_{p_2,r_2})}+
    \|p\|_{\widetilde{L}^\infty_T(B^{s_2}_{p_2,r_2}) }<\infty.
  \end{align*}
  Theorem \ref{a1} then enables us to extend the solution beyond $T.$ 
  This completes the proof of the theorem.\qed

 \section{Proof of Theorem \ref{a3}}
  First, applying Lemma \ref{dts} to the first equation of the $\widetilde{NSNPP}$ system
  \begin{align*}
  (u_{\nu})_t+u_\nu \cdot \nabla u_\nu-\nu\triangle u_\nu+\Pi(u_\nu,u_\nu)=\mathcal{P}\big((n_\nu-p_\nu)\psi_\nu\big)
  \end{align*}
  with $\rho=\infty,~\rho_1=1,$ we have
  \begin{align*}
  \|u_{\nu}\|_{\widetilde{L}^\infty_t(B^{s_1}_{p_1,r_1})}\lesssim & exp(C\int_0^t\|u_{\nu}\|_{B^{s_1}_{p_1,r_1}}dt')\Big(
 \|u_{\nu0}\|_{B^{s_1}_{p_1,r_1}}\\\nonumber&+\|\Pi(u_{\nu},u_{\nu})\|_{\widetilde{L}^1_t(B^{s_1}_{p_1,r_1})}+
  \|\mathcal{P}\big((n_\nu-p_\nu)\psi_\nu\big)\|_{\widetilde{L}^1_t(B^{s_1}_{p_1,r_1})}\Big).
\end{align*}
This works in the same way as applying Lemma \ref{ts} to the first equation of the $\widetilde{ENPP}$ system
\begin{align*}u_t+u\cdot \nabla u+\Pi(u,u)=\mathcal{P}
\big((n-p)\psi\big)\end{align*}
to get
\begin{align*}
  \|u\|_{\widetilde{L}^\infty_t(B^{s_1}_{p_1,r_1})}\lesssim & exp(C\int_0^t\|u\|_{B^{s_1}_{p_1,r_1}}dt')\Big(
 \|u_{0}\|_{B^{s_1}_{p_1,r_1}}\\\nonumber&+\|\Pi(u_{},u_{})\|_{\widetilde{L}^1_t(B^{s_1}_{p_1,r_1})}+
  \|\mathcal{P}\big((n-p)\psi\big)\|_{\widetilde{L}^1_t(B^{s_1}_{p_1,r_1})}\Big)
\end{align*}
in Section 3. Thus following along the same lines as above, we conclude that there exist positive constants $T$ and $M$ independent of $\nu,$ such that the $NSNPP$ system has a unique solution $(u_\nu,n_\nu,p_\mu,P_\nu,\Phi_\nu)$ in $X(T)\times Y_\alpha(T)~ $ for some $\alpha\in(0,1),$ and $\|(u_\nu,n_\nu,p_\nu)\|_{X(T)}\leq M.$

Next, let
$(u,n,p)$ and $(u_\nu,n_\nu,p_\mu)\in X(T)$
be the solutions of the $ENPP$ system and the $NSNPP$ system respectively with the same initial data. The difference $u_\nu-u$ satisfies
\begin{align*}
     (u_\nu-u)_t&-\nu\triangle(u_\nu-u)=(u_\nu-u)\nabla u_\nu-\Pi(u_\nu-u,u_\nu+u)\\&+\mathcal{P}\big((n_\nu-p_\nu)(\psi_\nu-\psi)\big)
   +\mathcal{P}\big((n_\nu-n-p_\nu+p)\psi)+u\cdot \nabla (u_\nu-u\big)+\nu\triangle u.
\end{align*}
Lemma \ref{dts} implies
\begin{align*}
  \|u_\nu-u\|_{\widetilde{L}^\infty_t(B^{s'_1}_{p_1,r_1})}\lesssim &
  exp(C\int_0^t\|u_{\nu}\|_{B^{s_1}_{p_1,r_1}}dt')
  \Big(\|\Pi(u_\nu-u,u_\nu+u)\|_{\widetilde{L}^1_t(B^{s'_1}_{p_1,r_1})}
    \\\nonumber&+\|\mathcal{P}\big((n_\nu-n-p_\nu+p)\psi\big)\|_{\widetilde{L}^1_t(B^{s'_1}_{p_1,r_1})}
    +\|\mathcal{P}\big((n_\nu-p_\nu)(\psi_\nu-\psi)\big)\|_{\widetilde{L}^1_t(B^{s'_1}_{p_1,r_1})}\\\nonumber&
    +(1+\nu t)^{\frac{1}{2}}\nu^{-\frac{1}{2}}\|\nu\triangle u\|_{\widetilde{L}^2_t(B^{s'_1-1}_{p_1,r_1})}\Big).
\end{align*}
It is easy to obtain that if $\nu\leq 1,$
\begin{align*}(1+\nu t)^{\frac{1}{2}}\nu^{-\frac{1}{2}}\|\nu\triangle u\|_{\widetilde{L}^2_t(B^{s'_1-1}_{p_1,r_1})}\lesssim (1+\nu T)^{\frac{1}{2}}\nu^{\frac{1}{2}}\| u\|_{\widetilde{L}^2_t(B^{s_1}_{p_1,r_1})}\lesssim (1+ T)^{\frac{1}{2}}T^{\frac{1}{2}}\nu^{\frac{1}{2}}\| u\|_{\widetilde{L}^\infty_T(B^{s_1}_{p_1,r_1})}.
\end{align*}
Then reasoning along exactly the same lines as that of Lemma \ref{mm}, we get
 \begin{align*}
    f(t)\lesssim C(T,M)\big(\nu^{\frac{1}{2}}+\int_0^tf(t')dt'\big),
  \end{align*}
  where $$f(t)\triangleq \|u_\nu-u\|_{\widetilde{L}^\infty_t(B^{s'_1}_{p_1,r_1})}+\|n_\nu-n \|_{\widetilde{L}^\infty_t(B^{s_2-1}_{p_2,r_2})}+
  \|p_\nu-p\|_{\widetilde{L}^\infty_t(B^{s_2-1}_{p_2,r_2})},$$ and $C(T,M)$ is a constant depending on $T$ and $M.$ Applying the Gronwall lemma completes the proof of Theorem \ref{a3}.\qed

\bigskip
\noindent\textbf{Acknowledgements}. This work was partially
supported by NNSFC (No. 11271382), RFDP (No. 20120171110014), and
the key project of Sun Yat-sen University.
\phantomsection
\addcontentsline{toc}{section}{\refname}
\bibliography{reference}

\end{document}